\documentclass[11pt,epsfig]{article}
\usepackage{amsmath,amssymb,amsfonts}
\usepackage{amsthm}
\usepackage{epsfig}
\usepackage{graphicx}
\usepackage{mathrsfs}
\usepackage{hyperref}
\usepackage{color}
\usepackage[margin=1.2in]{geometry}
\numberwithin{equation}{section}
\usepackage{cite}
\usepackage{subfigure}
\usepackage{lipsum,titletoc}
\titlecontents{section} % set formatting for \section -
[2.5em]                 % adjust left margin
{\rmfamily}%{\scshape}             % font formatting
{\contentslabel{2.3em}} % section label and offset
{\hspace*{-2.3em}} {\titlerule*[1pc]{.}\contentspage}

\titlecontents{subsection} % set formatting for \section -
[4.8em]                 % adjust left margin
{\rmfamily}             % font formatting
{\contentslabel{2.3em}} % section label and offset
{\hspace*{-2.3em}} {\titlerule*[1pc]{.}\contentspage}

\definecolor{darkred}{rgb}{0.7,0,0}

\definecolor{green}{rgb}{0,0.7,0}

\newtheoremstyle{thmm}{1.5ex plus 1ex minus .2ex}{1.5ex plus
1ex minus
.2ex}{\rmfamily}{}{\bfseries}{}{1em}{} \theoremstyle{thmm}
\newtheorem{theorem}{Theorem}[section]
\newtheorem{lemma}{Lemma}[section]

\newtheorem{definition}{Definition}[section]
\newtheorem{remark}{Remark}[section]
\renewenvironment{proof}[1][Proof]{\noindent\textit{#1. }
}{\hfill$\square$}

\allowdisplaybreaks
\renewcommand{\theequation}{\thesection.\arabic{equation}}

\newcommand{\nn}{\nonumber}
\def \endproof{\vrule height8pt width 5pt depth 0pt}

\newcommand{\vertiii}[1]
{{\left\vert\kern-0.25ex\left
\vert\kern-0.25ex\left\vert #1
    \right\vert\kern-0.25ex\right
\vert\kern-0.25ex\right\vert}}

\def\R{\mathbb{R}}
\def\C{\mathbb{C}}
\def\d{{\rm d}}

\begin{document}

\title{\bf%\Large 
Convergence of a decoupled mixed FEM for 
the dynamic Ginzburg--Landau equations 
in nonsmooth domains 
with incompatible initial data
\thanks{The research stay of the author at Universit\"at T\"ubingen 
was supported by the Alexander von Humboldt Foundation. 
This work was supported in part by the NSFC (grant no. 11301262).
}
} 
\author{ 
Buyang Li 
\footnote{
%Department of Applied Mathematics, 
%The Hong Kong Polytechnic University, Kowloon, Hong Kong. 
%}
%\footnote{
Mathematisches Institut, 
Universit\"at T\"ubingen, 
72076 T\"ubingen, Germany. 
{\tt li@na.uni-tuebingen.de} }
}
\date{}
\maketitle 
\vspace{-15pt}

\begin{abstract}
In this paper, we propose a fully discrete mixed finite element method for solving the time-dependent Ginzburg–Landau equations, and prove the convergence of the finite element solutions in general curved polyhedra, possibly nonconvex and multi-connected, without assumptions on the regularity of the solution.  
Global existence and uniqueness of weak solutions for 
the PDE problem are also obtained in the meantime. 
A decoupled time-stepping scheme is introduced, 
which guarantees that 
the discrete solution has bounded discrete energy,  
and the finite element spaces are chosen 
to be compatible with the nonlinear structure of the equations. 
Based on the boundedness of the discrete energy, 
we prove the convergence 
of the finite element solutions 
by utilizing a uniform $L^{3+\delta}$ regularity 
of the discrete harmonic vector fields, establishing 
a discrete Sobolev embedding inequality 
for the N\'ed\'elec finite element space, 
and introducing a $\ell^2(W^{1,3+\delta})$ estimate
for fully discrete solutions of parabolic equations. 
The numerical example shows that the constructed 
mixed finite element solution 
converges to the true solution of the PDE problem
in a nonsmooth and multi-connected domain, while the standard
Galerkin finite element solution does not converge. 
\end{abstract}

\pagestyle{myheadings}
\thispagestyle{plain}
%\markboth{}{}

\tableofcontents

\section{Introduction}
\setcounter{equation}{0} 

The time-dependent Ginzburg--Landau 
equation (TDGL) is a macroscopic phenomenological 
model for the superconductivity phenomena 
in both low and high temperatures 
\cite{GL,GE,Gennes,Tinkham},  
and has been widely accepted in the 
numerical simulation of transition and
vortex dynamics of both type-I and type-II superconductors  
\cite{FUD91,LMG91}.  
In a non-dimensionalization form, the TDGL is given by 
\begin{align}
&\eta\frac{\partial \psi}{\partial t} + i\eta\kappa\psi\phi
+ \bigg(\frac{i}{\kappa} \nabla + \mathbf{A}\bigg)^{2} \psi
 + (|\psi|^{2}-1) \psi  = 0,
\label{GLLPDEq1}\\[5pt]
&\frac{\partial \mathbf{A}}{\partial t}  +\nabla \phi
+ \nabla\times(\nabla\times{\bf A})
+  {\rm Re}\bigg[\overline\psi\bigg(\frac{i}{\kappa} \nabla 
+ \mathbf{A}\bigg) \psi\bigg] =  \nabla\times {\bf H} ,
\label{GLLPDEq2}
\end{align}  
where the order parameter $\psi$ is 
complex scalar-valued, 
the electric potential $\phi$ is real scalar-valued and 
magnetic potential ${\bf A}$ is real vector-valued;   
$\eta>0$ and $\kappa>0$ are physical parameters, and  
${\bf H}$ is a time-independent external magnetic field.
In a domain $\Omega\subset\R^3$
occupied by a superconductor, 
the following physical boundary conditions are 
often imposed:  
\begin{align} 
&\Big(\frac{i}{\kappa}\nabla\psi +
\mathbf{A}\psi\Big) \cdot{\bf n} = 0
&\mbox{on}\,\,\,\,\partial\Omega , \\
&{\bf n}\times{\bf B}={\bf n}\times {\bf H}
&\mbox{on}\,\,\,\,\partial\Omega , \\
&{\bf E}\cdot{\bf n}=0  &\mbox{on}\,\,\,\,\partial\Omega , 
\label{OrBDC3}
\end{align}
where $\mathbf{n}$ denotes the unit normal vector on 
the boundary of the domain, ${\bf B}=\nabla\times{\bf A}$ 
and ${\bf E}=-\partial_t{\bf A}-\nabla\phi$
denote the induced magnetic and electric fields, respectively.

Besides \eqref{GLLPDEq1}-\eqref{GLLPDEq2}, 
an additional gauge condition is needed 
for the uniqueness of the solution $(\psi,\phi,{\bf A})$.
Under the gauge $\phi=-\nabla\cdot{\bf A}$, 
the TDGL reduces to 
\begin{align}
&\eta\frac{\partial \psi}{\partial t} -i\eta \kappa \psi \nabla\cdot{\bf A}
+ \bigg(\frac{i}{\kappa} \nabla + \mathbf{A}\bigg)^{2} \psi
 + (|\psi|^{2}-1) \psi = 0,
\label{PDE1}\\[5pt]
&\frac{\partial \mathbf{A}}{\partial t} 
+ \nabla\times(\nabla\times{\bf A})
-\nabla(\nabla\cdot{\bf A}) 
+  {\rm Re}\bigg[\overline\psi\bigg(\frac{i}{\kappa} \nabla 
+ \mathbf{A}\bigg) \psi\bigg] =  \nabla\times {\bf H} ,
\label{PDE2}
\end{align}
and the boundary conditions can be written as 
\footnote{
Since \eqref{PDEBC} implies $\partial_t{\bf A}\cdot{\bf n}=0$, 
\eqref{PDEBC-00} and \eqref{PDEBC} imply   
${\rm Re}\big[\overline\psi\big(\frac{i}{\kappa} \nabla 
+ \mathbf{A}\big) \psi\big]\cdot{\bf n}=0$ 
and \eqref{PDEBC-0} implies 
$[\nabla\times(\nabla\times {\bf A}-{\bf H})]\cdot{\bf n}=0$ 
(if a vector field ${\bf u}$ satisfies ${\bf n} \times {\bf u} = 0$ 
on $\partial\Omega$, then
$(\nabla\times{\bf u}) \cdot{\bf n}= 0$ on $\partial\Omega$), 
it follows from \eqref{PDE2} that 
$\nabla\phi\cdot{\bf n}=-\nabla(\nabla\cdot{\bf A})\cdot{\bf n}=0$ on each 
smooth piece of $\partial\Omega$. 
Hence, \eqref{PDEBC-00}-\eqref{PDEBC} imply \eqref{OrBDC3}. }
\begin{align} 
&\nabla\psi\cdot{\bf n} = 0
&&\mbox{on}\,\,\,\,\partial\Omega ,  \label{PDEBC-00}\\
&{\bf n}\times(\nabla\times{\bf A})= {\bf n}\times{\bf H}  
&&\mbox{on}\,\,\,\,\partial\Omega , \label{PDEBC-0}\\
&{\bf A}\cdot{\bf n}=0 
&&\mbox{on}\,\,\,\,\partial\Omega .
\label{PDEBC}
\end{align}
Given the initial conditions
\begin{align}\label{PDEini}
\psi(x,0)=\psi_0(x)\quad\mbox{and}\quad
{\bf A}(x,0)={\bf A}_0(x),\quad\mbox{for}\,\, x\in\Omega , 
\end{align}
the solution $(\psi,{\bf A})$ can be solved 
from \eqref{PDE1}-\eqref{PDEini}. 
Other gauges can also be used, and
the solutions under different gauges 
are equivalent in the sense that they produce
the same quantities of physical intereset 
\cite{CHL,Tinkham}, such as 
the superconducting density $|\psi|^2$ and  
the magnetic field ${\bf B}$.

In a smooth domain,  
well-posedness of \eqref{PDE1}-\eqref{PDEini} has been proved
in \cite{CHL} and convergence of 
the Galerkin finite element method (FEM)
was proved in \cite{CD01,GLS} 
with different time discretizations 
by assuming that the PDE's solution is smooth enough,
e.g. ${\bf A}\in L^\infty(0,T;{\bf H}^1)\cap L^2(0,T;{\bf H}^2)$.
In a nonsmooth domain 
such as a curved polyhedron, 
the magnetic potential ${\bf A}$ may be only in 
$L^\infty(0,T;{\bf H}({\rm curl,div}))\cap 
L^2(0,T;{\bf H}^{1/2+\delta})$,
where $\delta>0$ can be arbitrarily small 
(depending on the angle of the edges
or corners of the domain), 
and so the Galerkin finite element solution 
may not converge to the solution of \eqref{PDE1}-\eqref{PDE2}. 
Some mixed FEMs 
were proposed in \cite{Chen97,GS15}, 
and the numerical simulations in \cite{GS15}
show better results in nonsmooth domains,  
compared with the Galerkin FEM. 
Some discrete gauge invariant numerical methods 
\cite{Du98,DJ05} are also promising
to approximate the solution correctly. 
Convergence of these numerical methods
have been proved in the case that the PDE's solution
is smooth enough. 
However, whether the numerical solutions converge to
the PDE's solution in nonsmooth domains where the magnetic 
potential is only in 
$L^\infty(0,T;{\bf H}({\rm curl,div}))\cap 
L^2(0,T;{\bf H}^{1/2+\delta})$ 
is still unknown. 
Another problem is that 
the initial data ${\bf A}_0$ are often incompatible with 
the boundary condition \eqref{PDEBC-0} 
(see the numerical examples in \cite{ASPM11,CD01,Mu97}, 
where ${\bf n}\times(\nabla\times{\bf A}_0)=0$ but ${\bf n}\times{\bf H}\neq 0$),
and this also leads to low regularity of the solution.

Numerical analysis of the TDGL 
under the zero electric potential gauge $\phi=0$ 
has also been done in many works 
\cite{ASPM11,GKL02,GKLLP96,MH98,
RPCA04,VMB03,WA02,Yang}; 
also see the review paper \cite{Du05}. 
Since 
$\|\nabla\times{\bf A}\|_{L^2 }$ is not equivalent to 
$\|\nabla{\bf A}\|_{L^2 }$, 
both theoretical and numerical analysis are difficult 
under this gauge without extra 
assumptions on the regularity of the 
PDE's solution. 
Again, convergence of these numerical methods
have been proved in the case that the PDE's solution
is smooth enough.

%However, convergence of the numerical solutions 
%has not been proved in nonsmooth domains 
%without extra regularity 
%assumptions on the PDE's solution. 
Under either gauge, 
convergence of the numerical solutions has not been
proved in nonsmooth domains such as general curved polyhedra, 
possibly nonconvex and multi-connected. 
Meanwhile, correct numerical approximations of the TDGL 
in domains with edges and corners are 
important for physicists 
and engineers \cite{ASPM11,BKP05,VMB03}.
The difficulty of the problem is to control
the nonlinear terms in the equations 
only based on the a priori estimates
of the finite element solution. 
In this paper, we introduce a decoupled
mixed FEM for solving \eqref{PDE1}-\eqref{PDEBC} 
which guarantees that 
the discrete solution has bounded discrete energy, 
and prove convergence of the fully discrete 
finite element solution in general curved polyhedra 
without assumptions on the regularity of the
PDE's solution. 
We control the nonlinear terms by proving 
a uniform $L^{3+\delta}$ regularity 
for the discrete harmonic vector
fields in curved polyhedra, 
establishing a discrete Sobolev compact 
embedding inequality ${\bf H}_h({\rm curl,div})\hookrightarrow
\hookrightarrow {\bf L}^{3+\delta}$
for the functions in the N\'ed\'elec element space,  
and introducing a $\ell^2(W^{1,3+\delta})$ estimate
for fully discrete finite element solutions of 
parabolic equations,
where $\delta>0$ is some constant which
depends on the given domain.

%\pagebreak 

\section{Main results}\label{MainR}
\setcounter{equation}{0}

\subsection{A decoupled mixed FEM
with bounded discrete energy}

In this subsection, we introduce our assumptions 
on the domain and define the fully discrete finite element method 
to be considered in this paper. 
Then we introduce a discrete energy function
(different from the free energy) 
and sketch a proof for a basic energy inequality 
satisfied by the finite element solution. \medskip  

\begin{figure}[htp]
\centering
%\subfigure[The domain $\Omega$ (shadow part)] 
{\includegraphics[height=1.3in,width=1.5in]{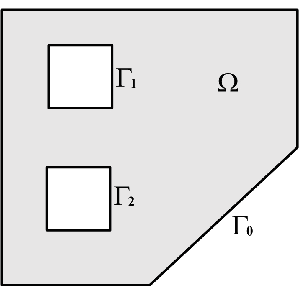}}
\qquad\qquad
%\subfigure[The domain $\Omega_0$ (shadow part)] 
{\includegraphics[height=1.3in,width=1.5in]{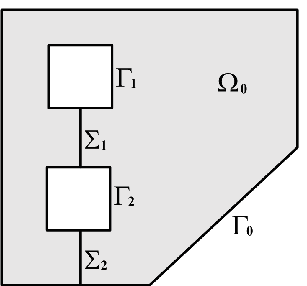}}
\caption{Illustration of the domain 
{\rm(}$\Omega$ is the shadow region{\rm)}.}
\label{FigD}
\end{figure}

\begin{definition}
A {\it curved polyhedron (or polygon)} 
is a bounded Lipschitz domain $\Omega\subset\R^3$ 
(or $\Omega\subset\R^2$), 
possibly nonconvex and multi-connected, 
such that its boundary is locally 
$C^\infty$-isomorphic to the boundary of a polyhedron \cite{BS87},  
and there are $\frak M$ 
pieces of surfaces $\Sigma_1$,
$\cdots$, $\Sigma_{\frak M}$ transversal
to $\partial\Omega$ such that $\Sigma_i\cap\Sigma_j=\emptyset$
for $i\neq j$ and the domain 
$\Omega_0:=\Omega\backslash\Sigma$
is simply connected, 
where $\Sigma=\cup_{j=1}^{\frak M}\Sigma_j$ 
(see Figure \ref{FigD}) . \medskip
\end{definition}

\begin{remark}
The integer $\frak M$ is often referred to as the
first Betti number of the domain. 
The existence of the surfaces 
$\Sigma_j$, $j=1,\cdots,\frak M$,  
is only needed in the analysis of the 
finite element solutions
by using the Hodge decomposition \cite{KY09}. 
One does not need to know these surfaces
in practical computation. 
\end{remark}

{\bf Assumptions 2.1.}$\,\,\,$
We assume that $\Omega\subset\R^3$ is a curved polyhedron 
which is partitioned into quasi-uniform tetrahedra. 
%which fit the boundary $\partial\Omega$ exactly.
For any given integers 
\begin{align}\label{Condrk}
r\geq 1 
\quad\mbox{and}\quad
k\geq 2r-1 ,
\end{align} 
we denote by ${\mathbb S}_h^{r}$ 
the complex-valued Lagrange finite 
element space of degree $\leq r$, 
denote by ${\mathbb V}_h^{k+1}$ 
the real-valued Lagrange finite 
element space of degree $\leq k+1$,  
and let ${\mathbb N}_h^k$ be either the 
N\'ed\'elec 1st-kind H(curl) element space of order $k$ 
\cite{Ned80} 
or the N\'ed\'elec 
2nd-kind H(curl) element space of degree $\leq k$ \cite{Ned86} 
(also see page 60 of \cite{AFW}). 
\endproof \medskip

Let the time interval $[0,T]$
be partitioned into 
$0=t_0<t_1<\cdots<t_N$ uniformly,
with $\tau=t_{n+1}-t_n$. 
For any given functions $f_n$, $n=0,1,\cdots,N$,
we define its discrete time derivative as 
\begin{align*}
D_\tau f^{n+1}:=(f^{n+1}-f^n)/\tau ,\quad
n=0,1,\cdots,N-1 .
\end{align*}
We introduce a decoupled 
backward Euler scheme for solving 
\eqref{PDE1}-\eqref{PDE2}:  
\begin{align} 
&\eta D_\tau \psi^{n+1} 
-i\eta\kappa\Theta(\psi^{n})\nabla\cdot{\bf A}^n
+ \bigg(\frac{i}{\kappa} \nabla + {\bf A}^{n+1}\bigg)^{2} \psi^{n+1}
 + (|\psi^{n+1}|^{2}-1) \psi^{n+1}= 0, 
\label{TD1}\\[5pt] 
& D_\tau {\bf A}^{n+1} 
-\nabla(\nabla\cdot{\bf A}^{n+1})
+ \nabla\times(\nabla\times{\bf A}^{n+1})
+  {\rm Re}\bigg[\overline\psi^{n}\bigg(\frac{i}{\kappa} \nabla 
+{\bf A}^{n}\bigg) \psi^{n}\bigg]
=  \nabla\times {\bf H}  , 
\label{TD2}  
\end{align} 
where
we have used a cut-off function 
\begin{align}
\Theta(z):=z/\max(|z|,1) ,\quad\forall~ z\in\C ,
\end{align}
which satisfies $\Theta(z)=z$ if $|z|\leq 1$. 
%Since the solution of  
%\eqref{PDE1}-\eqref{PDE2} satisfies
%$|\psi|\leq 1$ (see Appendix), 
%our time discretization is consistent
%with the PDE problem. 

For any given integers $r$ and $k$
which satisfy the condition \eqref{Condrk},
we solve \eqref{TD1} by the Galerkin FEM 
and solve \eqref{TD2} by a mixed FEM.
Let $(\psi_{h}^{0}, {\bf A}_h^0):=(\psi_{0},
{\bf A}_0)$ at the initial time step and define 
$\phi_h^0:=-\nabla\cdot{\bf A}_0$. 
%be the 
%solution of the finite element equation 
%\begin{align} 
%& (\phi_h ^0, \chi_h) 
%- ({\bf A}_h^0,\nabla \chi_h ) = 0 \, , 
%\qquad\forall\, \chi_h\in {\mathbb V}_h^{k+1} .
%\label{FEM3} 
%\end{align} 
We look for 
$\psi_{h}^{n+1} \in {\mathbb S}_{h}^r$
and $(\phi_{h}^{n+1},{\bf A}_h ^{n+1} ) 
\in {\mathbb V}_{h}^{k+1} \times {\mathbb N}_{h}^{k}$,
$n=0,1,\cdots,N-1$, 
satisfying the equations 
\begin{align} 
&(\eta D_{\tau}\psi_{h}^{n+1}, \varphi_h)  
+ \bigg( \bigg(\frac{i}{\kappa}\nabla+{\bf A}_h^{n+1}\bigg) {\psi}_{h}^{n+1} \,, 
\bigg(\frac{i}{\kappa}\nabla+{\bf A}_h^{n+1}\bigg) \varphi_h\bigg) \nn \\
&\quad 
+((|\psi_{h}^{n+1}|^{2}-1)\psi_{h}^{n+1}, \varphi_h) 
=- (i\eta \kappa\Theta(\psi_{h}^{n})\phi_{h}^{n} ,\varphi_h) \, , 
&&\forall\, \varphi_h \in {\mathbb S}_{h}^r ,
\label{FEM1} \\[10pt]
& (\phi_h ^{n+1}, \chi_h) 
- ({\bf A}_h^{n+1},\nabla \chi_h ) = 0 \, , 
&&\forall\, \chi_h\in {\mathbb V}_{h}^{k+1}, 
\label{FEM3} \\[10pt]
&(D_{\tau}{\bf A}_{h}^{n+1},{\bf a}_{h}) 
+(\nabla\phi_{h}^{n+1} \, , {\bf a}_{h} )
+ (\nabla\times{\bf A}_h ^{n+1} \, , \nabla\times {\bf a}_{h}) \nn \\
& 
= ({\bf H} \, , \nabla\times {\bf a}_{h}) 
-{\rm Re}\bigg(\overline\psi_h^{n}\bigg(\frac{i}{\kappa} \nabla 
+{\bf A}_h^{n}\bigg) \psi_h^{n} , {\bf a}_h\bigg) \, , 
  &&\forall\, {\bf a}_h\in  {\mathbb N}_{h}^{k} .
\label{FEM2} 
\end{align}
After solving $\psi_h^{n+1}$, $\phi_h^{n+1}$
and ${\bf A}_h^{n+1}$ from the equations above,
the magnetic and electric fields
can be computed by 
${\bf B}_h^{n+1}=\nabla\times{\bf A}_h^{n+1} $
%\qquad\mbox{and}\qquad 
and 
${\bf E}_h^{n+1}=-D_\tau{\bf A}_h^{n+1}
-\nabla\phi_h^{n+1} $.  

\begin{remark}
For simplicity, we have chosen 
$(\psi_h^0,{\bf A}_h^0)=(\psi_0,{\bf A}_0)$
at the initial step, which are not finite element functions. 
Due to the nonlinearities and the choice of the initial data,   
some integrals in \eqref{FEM1} and \eqref{FEM2} may 
need to be evaluated numerically in practical computations. 
In this paper, we focus on the analysis 
of the discretization errors of the finite element method 
and assume that all the integrals 
are evaluated accurately. %\medskip
\end{remark}

\begin{remark}
Since we have not assumed any extra regularity 
of the PDE's solution, we need the condition \eqref{Condrk} 
to be compatible with the nonlinear structure
of the equations in order to control a nonlinear term
arising from \eqref{FEM1} 
(see \eqref{CondrkN} for the details).
If the PDE's solution is smooth enough, 
(e.g. consider the problem in a smooth domain), 
then the condition 
\eqref{Condrk} can be relaxed.  
\end{remark}

%\subsection{The discrete energy }
We define the discrete energy 
\begin{align}
{\cal G}_{h}^{n} 
&=\int_\Omega 
\bigg(\frac{1}{2}\bigg|\frac{i}{\kappa} 
\nabla\psi_h^{n}  + {\bf A}_h^{n}\psi_h^{n} \bigg|^2
+\frac{1}{4}\big(|\psi_h^{n}|^2-1)^2 \bigg)\d x \\
&\quad +
\int_\Omega
\bigg(\frac{1}{2} |\nabla\times {\bf A}_h^{n}-{\bf H}|^2
+\frac{1}{2} |\phi_h^{n}|^2 \bigg)\d x \nn 
\end{align}
for $n=0,1,\cdots,N$.  
By substituting $\varphi_h=D_\tau\psi_h^{n+1}$,
$\chi_h=\phi_h ^{n+1}$
and ${\bf a}_h=D_\tau{\bf A}_h^{n+1}$ 
into \eqref{FEM1}-\eqref{FEM2}, we obtain
\begin{align}\label{FEMEGDecay}
&D_\tau {\cal G}_h^{n+1}
+\int_\Omega 
\big((\eta-\tau/2) |D_\tau\psi_h^{n+1}|^2
+|D_\tau{\bf A}_h^{n+1}|^2 
 \big) \d x  \nn\\
&\quad 
+\int_\Omega 
\frac{\tau}{2}\bigg( \bigg|D_\tau\bigg(\frac{i}{\kappa} 
\nabla\psi_h^{n+1} + {\bf A}_h^{n}\psi_h^{n+1} \bigg)\bigg|^2
+\big|D_\tau\phi_h^{n+1}|^2 
+\big|D_\tau\nabla\times {\bf A}_h^{n+1}|^2 \bigg)\d x \nn\\
&\quad 
+\int_\Omega \frac{\tau}{2}\bigg(|\psi_h^{n+1}|^2|D_\tau \psi_h^{n+1}|^2 
+\frac{1}{2}|D_\tau |\psi_h^{n+1}|^2|^2 \bigg)\d x \nn\\
&= 
-\int_\Omega
i\eta\kappa\Theta(\psi_h^{n})\phi_h^{n}D_\tau\psi_h^{n+1}\d x 
+{\rm Re}\int_\Omega \tau D_\tau\bigg(\frac{i}{\kappa} 
\nabla\psi_h^{n+1} + {\bf A}^{n+1}\psi_h^{n+1}\bigg) 
\overline\psi_h^{n}D_\tau{\bf A}_h^{n+1} \d x \nn\\
&\leq \frac{\eta}{2}\int_\Omega|D_\tau\psi_h^{n+1}|^2\d x
+ \frac{\eta\kappa^2}{2}\int_\Omega 
|\phi_h^{n}|^2\d x  \nn\\
&\quad 
+\int_\Omega \frac{\tau}{2} \bigg|D_\tau\bigg(\frac{i}{\kappa} 
\nabla\psi_h^{n+1} + {\bf A}_h^{n}\psi^{n+1}\bigg)\bigg|^2\d x 
+\int_\Omega \frac{\tau}{2}|\psi_h^{n}|^2|D_\tau{\bf A}_h^{n+1}|^2 \d x , 
%\nn\\
%&\leq \frac{\eta}{2}\int_\Omega|D_\tau\psi_h^{n+1}|^2\d x
%+ \eta\kappa^2{\cal G}_h^{n} \nn\\
%&\quad +\int_\Omega \frac{\tau}{2} \bigg|D_\tau\bigg(\frac{i}{\kappa} 
%\nabla\psi_h^{n+1} + {\bf A}_h^{n}\psi_h^{n+1}\bigg)\bigg|^2\d x 
%+\int_\Omega \frac{\tau}{2}|\psi_h^{n}|^2|D_\tau{\bf A}_h^{n}|^2 \d x,
\end{align} 
which reduces to
\begin{equation}\label{FEMEGDecay0}
\begin{aligned}
&D_\tau {\cal G}_h^{n+1}
+\int_\Omega 
\bigg(\frac{\eta-\tau}{2} |D_\tau\psi_h^{n+1}|^2
+\frac{1}{2}|D_\tau{\bf A}_h^{n+1}|^2
+\frac{1-\tau|\psi_h^n|^2}{2}|D_\tau{\bf A}_h^{n+1}|^2
 \bigg) \d x  \\
 &\leq \eta\kappa^2{\cal G}_h^{n}  .
\end{aligned} 
\end{equation} 
Unlike the PDE's solution,
it is not obvious whether 
the finite element solution satisfies 
$|\psi_h^n|\leq 1$ pointwisely. 
In Section \ref{SecUnif} we shall prove 
\begin{align}\label{BDPsiDecay0}
1-\tau|\psi_h^n|^2\geq 0  
\end{align} 
when $\tau<\tau_0$ 
(for some positive constant $\tau_0$
which is independent of $\tau$ and $h$). 
Then \eqref{FEMEGDecay0} implies 
boundedness of the discrete energy 
via the discrete Gronwall's inequality. 
%reduces to
%\begin{align}\label{FEMEGDecay00}
%&D_\tau {\cal G}_h^{n+1}
%+\int_\Omega 
%\bigg(\frac{\eta-\tau}{2} |D_\tau\psi_h^{n+1}|^2
%+\frac{1}{2}|D_\tau{\bf A}_h^{n+1}|^2
% \bigg) \d x  \leq \eta\kappa^2{\cal G}_h^{n} .
%\end{align} 
%By using the discrete Gronwall's inequality,
%the last inequality implies that 
%the discrete energy  remains bounded
%as $\tau\rightarrow 0$, i.e. 
%\begin{align}\label{FEMEGDecay000}
%&\max_{0\leq n\leq N-1}{\cal G}_h^{n+1}
%+\sum_{n=0}^{N-1}\int_\Omega 
%\bigg(\frac{\eta-\tau}{2} |D_\tau\psi_h^{n+1}|^2
%+\frac{1}{2}|D_\tau{\bf A}_h^{n+1}|^2
% \bigg) \d x  \leq C{\cal G}_h^{0}  .
%\end{align} 
By utilizing the discrete energy, 
we derive further estimates  
which are used to prove compactness and convergence
of the finite element solution.

%Once  is done is Section \ref{SecUnif}.  
%Then the inequality above implies 
%boundedness of the discrete energy 
%(via the discrete Gronwall's inequality) 
%\begin{align}\label{FEEnergyG}
%&\max_{0\leq n\leq N-1}{\cal G}_h^{n+1} 
%+\sum_{n=0}^N\tau \int_\Omega 
%\bigg(\frac{\eta-\tau}{2} |D_\tau\psi_h^{n+1}|^2
%+\frac{1}{2}|D_\tau{\bf A}_h^{n+1}|^2 \bigg) \d x
%\leq C_T{\cal G}_h^0 
%%\leq e^{\eta\kappa^2T}{\cal G}^{**}(\psi_h^0,{\bf A}_h^0,{\bf B}_h^0)
%\end{align} 
%when $\tau<\min(\tau_0,\eta/2)$. 
%Meanwhile, we prove convergence of the 
%fully discrete finite element solution
%defined by \eqref{FEM1}-\eqref{FEM2} 
%by utilizing \eqref{FEEnergyG}. 

\subsection{Main theorem}

Let $W^{s,p}$, $s\geq 0$ and $1\leq p\leq\infty$, 
be the conventional Sobolev spaces of real-valued functions
defined on $\Omega$,
and let ${\bf W}^{s,p}=W^{s,p}\times W^{s,p}\times W^{s,p}$
be the corresponding Sobolev space of vector fields.  
The case of integer $s$ can be found in \cite{Adams},  
and the characterization of more
general function spaces 
with fractional $s$ can be found in \cite{Rychkov}. 
Let ${\cal W}^{s,p}:=W^{s,p}+i\,W^{s,p}$ denote the complex-valued Sobolev 
space and define the abbreviations 
\begin{align*}
&L^p:=W^{0,p},\quad\,  
{\cal L}^p:={\cal W}^{0,p},\quad\, 
{\bf L}^p:={\bf W}^{0,p} ,
\quad\mbox{for $1\leq p\leq\infty$,}\\
&H^s:=W^{s,2},\quad 
{\cal H}^s:={\cal W}^{s,2},\quad 
{\bf H}^s:={\bf W}^{s,2},
\quad\mbox{for $s\geq 0$} . 
\end{align*}
Moreover, we define
\begin{align}
&{\bf H}({\rm curl}):=
\{{\bf u}\in {\bf L}^2:\,
\nabla\times {\bf u}\in {\bf L}^2\}  , \\
& {\bf H}({\rm div})\,\,:=
\{{\bf u}\in {\bf L}^2:\,
\nabla\cdot {\bf u}\in L^2 \}  ,\\
%&{\bf H}({\rm curl,div}):=
%\{{\bf u}\in {\bf L}^2:\,
%\nabla\times {\bf u}\in {\bf L}^2\,\,\,
%\mbox{and}\,\,\,\nabla\cdot {\bf u}\in L^2\} ,\\
&{\bf H}({\rm curl,div}):=
\{{\bf u}\in {\bf L}^2:\,
\nabla\times {\bf u}\in {\bf L}^2,
\,\,\,\nabla\cdot {\bf u}\in L^2 \,\,\,\mbox{and}\,\,\, 
{\bf u}\cdot{\bf n}=0\,\,\mbox{on}\,\,\partial\Omega\} . \label{Def_Hcd_bd1}
%\\
%&{\bf H}^0({\rm curl,div}):=
%\{{\bf u}\in {\bf L}^2:\,
%\nabla\times {\bf u}\in {\bf L}^2\,\,\,
%\mbox{and}\,\,\,\nabla\cdot {\bf u}\in L^2,\,\,\,
%{\bf u}\times{\bf n}=0\,\,\mbox{on}\,\,\partial\Omega\} . 
%\label{Def_Hcd_bd2}
\end{align}

Let $\psi_{h,\tau}^+$, 
$\psi_{h,\tau}^-$, ${\bf A}_{h,\tau}^+$, 
${\bf A}_{h,\tau}^-$,
${\bf B}_{h,\tau}^+$ 
and ${\bf E}_{h,\tau}^+$ 
be the piecewise constant functions on $(0,T]$
such that on each subinterval $(t_n,t_{n+1}]$ 
\begin{align}
&\psi_{h,\tau}^+(t)=\psi_h^{n+1},  
&&\psi_{h,\tau}^-(t)=\psi_h^{n},  \\
&{\bf A}_{h,\tau}^+(t)={\bf A}_h^{n+1},   
&&{\bf A}_{h,\tau}^-(t)={\bf A}_h^{n} ,  \\
&{\bf B}_{h,\tau}^+(t)={\bf B}_h^{n+1}:=
\nabla\times{\bf A}_{h}^{n+1},
&&{\bf E}_{h,\tau}^+(t)={\bf E}_h^{n+1}
:=-D_\tau{\bf A}_h^{n+1}-\nabla\phi_h^{n+1}  .
\label{Eht1} 
\end{align}

In this paper we prove the following theorem. \medskip

\begin{theorem}\label{MainTHM}
{\it Under Assumption {\rm 2.1}, 
for any given $\psi_0\in {\cal H}^1$ 
and ${\bf A}_0\in {\bf H}({\rm curl,div})$
such that $|\psi_0|\leq 1$, 
the system 
\eqref{FEM1}-\eqref{FEM2} 
has a unique finite element solution 
when $\tau<\eta$ $(\eta$ is the parameter in \eqref{GLLPDEq1}{\rm )}, which 
converges to the unique solution
of \eqref{PDE1}-\eqref{PDEini}  
as $\tau,h\rightarrow 0$ 
in the following sense:
\begin{align}
&\psi_{h,\tau}^+\rightarrow \psi 
&&\mbox{strongly in\, $L^\infty(0,T;{\cal L}^2)$
and weakly$^*$\, in\, $L^\infty(0,T;{\cal H}^1)$ }, \\ 
%& 
%&&\mbox{weakly\,\,\, in\, $L^2(0,T;L^\infty)$} ,\\
%& 
%&&\mbox{weakly$^*$\, in\, $L^\infty(0,T;{\cal H}^1)$ } ,\\[5pt] 
&{\bf A}_{h,\tau}^+\rightarrow {\bf A}
&& \mbox{strongly in\, $L^\infty(0,T;{\bf L}^2)$ 
and weakly$^*$ in\, $L^\infty(0,T;{\bf H}({\rm curl}))$} , \\
%& 
%&& \mbox{weakly$^*$ in\, $L^\infty(0,T;{\bf H}({\rm curl}))$ } , \\
&\phi_{h,\tau}^+\rightarrow \phi 
%=-\nabla\cdot{\bf A}
&&\mbox{strongly\,\,\, in\, $L^2(0,T;L^2)$ 
and weakly$^*$ in\, $L^\infty(0,T;L^2)$}, \\ 
%& 
%&&\mbox{weakly$^*$ in\, $L^\infty(0,T;L^2)$} , \\
&{\bf B}_{h,\tau}^+\rightarrow {\bf B}
&&\mbox{weakly$^*$\, in\, $L^\infty(0,T;{\bf L}^2)$},\\
&{\bf E}_{h,\tau}^+\rightarrow {\bf E}
&&\mbox{weakly\,\,\, in\, $L^2(0,T;{\bf L}^2)$} . 
\end{align}
}
\end{theorem}

In the meantime of proving Theorem \ref{MainTHM},
we also obtain global well-posedness of the 
PDE problem \eqref{PDE1}-\eqref{PDEini} 
(see Appendix).

\begin{remark}
If $\Omega$ is a curved polygon in $\R^2$
and the external magnetic field ${\bf H}$ is perpendicular
to the domain, i.e. ${\bf H}=(0,0,H)$, then 
\eqref{PDE1}-\eqref{PDE2} hold when ${\bf H}$ is 
replaced by $H$, with 
the following two-dimensional notations: 
\begin{align*} 
&\nabla\times {\bf A}
=\frac{\partial A_2}{\partial x_1}-\frac{\partial A_1}{\partial x_2},
\qquad 
\nabla\cdot {\bf A}=\frac{\partial A_1}{\partial x_1}
+\frac{\partial A_2}{\partial x_2},\\
&\nabla\times H=\bigg(\frac{\partial H}{\partial x_2},\,
-\frac{\partial H}{\partial x_1}\bigg),\quad 
\nabla\psi=\bigg(\frac{\partial \psi}{\partial x_1},\,
\frac{\partial \psi}{\partial x_2}\bigg).
\end{align*}
With these notations, 
\eqref{FEM1}-\eqref{FEM2} can also be used for 
solving the two-dimensional problem,
and Theorem \ref{MainTHM} can also be proved 
in the similar way. 
\end{remark}

\subsection{An overview of the proof}\label{OverviewP}
%\setcounter{equation}{0}
%
%For the convenience of the readers,
%in this section, we present an overview
%of our proof for the convergence
%of the finite element solution in arbitrary given
%curved polyhedra. 
%Readers can also pass this section and
%read the details in Section 3--4 directly. 

Our basic idea is
%, similar as \cite{BLP09}, 
to introduce $\psi_{h,\tau}$ and ${\bf A}_{h,\tau}$
($\psi_{h,\tau}^+$ and ${\bf A}_{h,\tau}^+$) 
as the piecewise linear (piecewise constant) interpolation
of the finite element solutions $\psi_h^{n+1}$ and ${\bf A}_h^{n+1}$
in the time direction, respectively, 
and denote by $\psi_{h,\tau}^-$ and ${\bf A}_{h,\tau}^-$ 
the piecewise constant interpolation
of $\psi_h^{n}$ and ${\bf A}_h^{n}$, respectively. 
Rewrite the finite element equations 
%defined on discrete time levels
as the following equations defined continuously in time:
\begin{align}
&\int_0^T\bigg[(\eta \partial_t\psi_{h,\tau}, \varphi_{h,\tau}) 
+\bigg( \bigg(\frac{i}{\kappa}\nabla+{\bf A}_{h,\tau}^+\bigg)\psi_{h,\tau}^+ \,, 
\bigg(\frac{i}{\kappa}\nabla+{\bf A}_{h,\tau}^+\bigg)\varphi_{h,\tau}\bigg)
 \bigg]\d t \nn\\
&\quad 
+((|\psi_{h,\tau}^+|^{2}-1)\psi_{h,\tau}^+,\varphi_{h,\tau})\bigg]\d t=
\int_0^T(i\eta \kappa
\Theta(\psi_{h,\tau}^-)\phi_{h,\tau}^- ,\varphi_{h,\tau})\d t   . \nn
\end{align}
\begin{align*}%\label{FFFFF1}
&\int_0^T\bigg[ (\phi_{h,\tau}^+, \chi_{h,\tau}) 
- ({\bf A}_{h,\tau}^+,\nabla \chi_{h,\tau} ) \bigg]\d t= 0 \, , \\[10pt]
&\int_0^T\bigg[ (\partial_t{\bf A}_{h,\tau},{\bf a}_{h,\tau}) 
+ (\nabla\phi_{h,\tau}^+ \, ,{\bf a}_{h,\tau})
+ (\nabla\times{\bf A}_{h,\tau}^+  \, , \nabla\times {\bf a}_{h,\tau})
\bigg]\d t\nn\\
&\quad +\int_0^T\bigg[ {\rm Re}\bigg(\overline\psi_{h,\tau}^- \bigg(\frac{i}{\kappa} \nabla 
+{\bf A}_{h,\tau}^-\bigg) \psi_{h,\tau}^-  , {\bf a}_{h,\tau}\bigg)\bigg]\d t
=\int_0^T\bigg[ (\nabla\times{\bf H} \, ,{\bf a}_{h,\tau})\bigg]\d t  \, . 
%\label{FFFFF2}
\end{align*}
If we can prove compactness and convergence of 
a subsequence of 
$\partial_t\psi_{h,\tau}$, $\psi_{h,\tau}^\pm$, 
$\partial_t{\bf A}_{h,\tau}$, ${\bf A}_{h,\tau}^\pm$,
$\phi_{h,\tau}^\pm$,  
and prove that the limits of any subsequence  
coincide with the PDE's solution, then we can conclude that 
the sequences $\psi_{h,\tau}^+$, 
and ${\bf A}_{h,\tau}^+$ converge to the PDE's solution 
as $h,\tau\rightarrow 0$.

To estimate the finite element solution
(in order to prove the compactness),
we introduce a discrete energy function ${\cal G}_{h}^{n} $
%\begin{align}
%{\cal G}_{h}^{n} 
%&=\int_\Omega 
%\bigg(\frac{1-\tau}{2}\bigg|\frac{i}{\kappa} 
%\nabla\psi_h^{n}  + {\bf A}_h^{n-1}\psi_h^{n} \bigg|^2
%+\frac{1}{4}\big(|\psi_h^{n}|^2-1)^2 \bigg)\d x \\
%&\quad +
%\int_\Omega
%\bigg(\frac{1}{2} |\nabla\times {\bf A}_h^{n}-{\bf H}|^2
%+\frac{1}{2} |\phi_h^{n}|^2 \bigg)\d x \nn\\
%&\quad +\int_\Omega \frac{\tau}{2}\bigg(\bigg|\bigg(\frac{i}{\kappa} 
%\nabla\psi_h^{n}+ {\bf A}_h^{n-1} \psi_h^{n}\bigg)
%+\psi_h^{n}D_\tau{\bf A}_h^{n}\bigg|^2
%+(1-|\psi_h^{n}|^2)|D_\tau{\bf A}_h^{n}|^2 \bigg)\d x   \nn
%\end{align}
%for $n=0,1,\cdots,N$,  
%with ${\bf A}_h^{-1}={\bf A}_h^0:={\bf A}_0$
%and $\phi_h^0:=-\nabla\cdot{\bf A}_0$,  
and a special time-stepping scheme 
from which one can derive \eqref{FEMEGDecay0}.   
%\begin{align}\label{FEMEGDecay01}
%&D_\tau {\cal G}_h^{n+1}
%+\int_\Omega 
%\bigg(\frac{\eta-\tau}{2} |D_\tau\psi_h^{n+1}|^2
%+\frac{1}{2}|D_\tau{\bf A}_h^{n+1}|^2
%+\frac{1-\tau|\psi_h^n|^2}{2}|D_\tau{\bf A}_h^{n}|^2
% \bigg) \d x  \leq \eta\kappa^2{\cal G}_h^{n} .
%\end{align} 
By proving \eqref{BDPsiDecay0}, 
we derive boundedness of the discrete energy
from \eqref{FEMEGDecay0} (via the discrete Gronwall's inequality). 
Based on the boundedness of the discrete energy, 
some further estimates need to be 
derived in order to prove convergence of the finite 
element solution. 
For example, 
in order to prove the weak convergence of 
a subsequence of 
$$
{\rm Re}\bigg[\overline\psi_{h,\tau}^- \bigg(\frac{i}{\kappa} \nabla 
+{\bf A}_{h,\tau}^-\bigg) \psi_{h,\tau}^-\bigg] 
$$
in $L^2(0,T;L^2)$, 
we need to prove the following convergence
(for a subsequence): 
\begin{align}
&\mbox{$\psi_{h,\tau}^-$ converges 
weakly in $L^2(0,T;W^{1,3+\delta})\hookrightarrow L^2(0,T;L^\infty)$
for some $\delta>0$}, \label{psihnL3}\\
&\mbox{$\psi_{h,\tau}^-$\,\, converges strongly in $L^\infty(0,T;L^{6-\epsilon})$ for arbitrarily small $\epsilon>0$}\\
&\mbox{${\bf A}_{h,\tau}^-$\, converges strongly in $L^\infty(0,T;L^{3+\delta})$
for some $\delta>0$} . \label{AAhnL3}
\end{align}
The boundedness of the discrete energy only 
implies the boundedness of  
$$
\|\psi_{h,\tau}^-\|_{L^\infty(0,T;H^1)} , \quad
\|{\bf A}_{h,\tau}^-\|_{L^\infty(0,T;L^2)}, \quad
\|\nabla\times {\bf A}_{h,\tau}^-\|_{L^\infty(0,T;L^2)} 
\quad
\mbox{and}\quad 
\|\phi_{h,\tau}^-\|_{L^\infty(0,T;L^2)} ,
$$ 
which are not enough for
$\psi_{h,\tau}^-$ and ${\bf A}_{h,\tau}^-$
to be compact and converge in the sense of \eqref{psihnL3}-\eqref{AAhnL3}. 

We shall prove \eqref{AAhnL3} by establishing 
a discrete Sobolev embedding inequality (Lemma 
\ref{SobolevD}):
\begin{align}\label{AhL33}
\|{\bf A}_{h}^n\|_{L^{3+\delta}}
\leq C(\|{\bf A}_{h}^n\|_{L^2}
+\|\nabla\times {\bf A}_{h}^n\|_{L^2}+\|\phi_{h}^n\|_{L^2}) ,
\end{align} 
and we also need to show that this embedding 
is compact.
Since we allow the domain to be multi-connected,
in order to prove \eqref{AhL33},
we need to use the discrete Hodge
decomposition
\begin{align*}
{\bf A}_h^n =
{\bf c}_h + \nabla\theta_h 
+ \sum_{j=1}^{\frak M}\alpha_{j,h}{\bf w}_{j,h}  
\end{align*}
and show that the divergence-free part ${\bf c}_h$,
the curl-free part $\nabla\theta_h $ and 
the discrete harmonic 
part $ \sum_{j=1}^{\frak M}\alpha_{j,h}{\bf w}_{j,h}$
are all bounded in ${\bf L}^{3+\delta}$. 
For this purpose, we need to construct 
the basis functions ${\bf w}_{j,h}$, $j=1,\cdots,\frak M$, 
of the discrete harmonic vector fields 
and prove that they are 
bounded in ${\bf L}^{3+\delta}$
(Lemma \ref{RegDHarV}).  

In order to prove \eqref{psihnL3},
we rewrite the finite element equation
of $\psi_h^{n+1}$ in the form of  
\begin{align*}%\label{EQDPSIa}
\eta D_\tau\psi_h^{n+1}
-\frac{1}{\kappa^2}\Delta_h\psi_h^{n+1}
=f_h^{n+1}
\end{align*}
%where 
%\begin{align*} 
%f_h^{n+1}=&
%-\frac{i}{\kappa}P_h\big(\nabla\psi^{n+1}_h\cdot{\bf A}_h^n\big)
%-\frac{i}{\kappa}\nabla_h\cdot\big(\psi^{n+1}_h {\bf A}_h^n\big)  &
%  \nn\\
%&-P_h\Big( |\mathbf{A}^{n}_h|^2\psi^{n+1}_h + 
%(|\psi^{n+1}_h|^{2}-1) \psi^{n+1}_h 
%+i\eta\kappa \Theta(\psi_h^{n})\phi_h^n\Big) ,
%\end{align*}
and prove the following inequality (Lemma 3.8): 
%in any given curved polyhedron there holds 
\begin{align}\label{LemW1qDa}
\sum_{n=0}^{N-1}\tau\|\psi_h^{n+1}\|_{W^{1,q+\delta}}^2\leq 
C\sum_{n=0}^{N-1}\tau\|f_h^{n+1}\|_{L^{q/2}}^2
+C\|\psi_h^0\|_{H^1}^2
\quad\mbox{ 
for some $q>3$ and $\delta>0$. }
\end{align} 
Then we prove 
\begin{align}
\sum_{n=0}^{N-1}\tau\|f_h^{n+1}\|_{L^{q/2}}^2
&\leq C+C\sum_{n=0}^{N-1}\tau\|\psi_h^{n+1}\|_{W^{1,q}}^2 \nn \\
&\leq C+C_\epsilon\sum_{n=0}^{N-1}\tau\|\psi_h^{n+1}\|_{H^{1}}^2
+\epsilon\sum_{n=0}^{N-1}\tau\|\psi_h^{n+1}\|_{W^{1,q+\delta}}^2  ,
\quad\forall\,\epsilon\in(0,1) .
\end{align} 
The last two inequalities imply
\begin{align}
\sum_{n=0}^{N-1}\tau\|\psi_h^{n+1}\|_{W^{1,q+\delta}}^2\leq 
C+C\sum_{n=0}^{N-1}\tau\|\psi_h^{n+1}\|_{H^{1}}^2
+C\|\psi_h^0\|_{H^1}^2
\leq C. 
\end{align} 

The compactness and convergence of the finite element 
solution are proved based on the uniform estimates established.
On one hand, in both \eqref{AhL33} and \eqref{LemW1qDa}  
we need some constant $\delta>0$
(which depends on the given curved polyhedron) to prove the 
convergence of the finite element solution. 
On the other hand, both \eqref{AhL33} and \eqref{LemW1qDa} are sharp: 
for any $\delta>0$ there exists a polyhedron such that
\eqref{AhL33} and \eqref{LemW1qDa} do not hold.

%Then we prove that the limit of
%$(\psi_{h,\tau},{\bf A}_{h,\tau})$, denoted by $(\Psi,{\bf\Lambda})$,  
%is a weak solution of the PDE problem
%and satisfies $|\Psi|\leq 1$ a.e. in $\Omega\times(0,T)$. 
%With this pointwise estimate and the regularity 
%established for $(\Psi,{\bf\Lambda})$, 
%the PDE problem has a unique weak solution,
%which implies   
%$(\Psi,{\bf\Lambda})=(\psi,{\bf A})$. 

\section{Proof of Theorem \ref{MainTHM}}\label{SecProof}
\setcounter{equation}{0}  

By substituting $\chi_h=\phi_h^{n+1}$ and 
${\bf a}_h={\bf A}_h^{n+1}$ into the equations
\begin{align} 
& (\phi_h ^{n+1}, \chi_h) 
- ({\bf A}_h^{n+1},\nabla \chi_h ) = 0 \, , 
&&\forall\, \chi_h\in {\mathbb V}_{h}^{k+1}, \\[10pt]
&\frac{1}{\tau} ({\bf A}_{h}^{n+1},{\bf a}_{h}) 
+(\nabla\phi_{h}^{n+1} \, , {\bf a}_{h} )
+ (\nabla\times{\bf A}_h ^{n+1} \, , \nabla\times {\bf a}_{h}) =0 ,
  &&\forall\, {\bf a}_h\in  {\mathbb N}_{h}^{k} ,
\end{align}
we see that the two equations above have only zero solution. 
Hence, 
for any given $(\psi_h^{n},{\bf A}_h^{n})\in {\mathbb S}_{h}^{r} \times {\mathbb N}_{h}^{k}$,  
the linear system \eqref{FEM3}-\eqref{FEM2} has a unique solution 
$(\phi_h^{n+1},{\bf A}_h^{n+1})\in 
{\mathbb V}_{h}^{k+1} \times {\mathbb N}_{h}^{k}$. 

Under the condition $\tau<\eta$, 
it is easy to see that for any given ${\bf A}_h^{n+1}\in{\mathbb N}_h^k$ the nonlinear operator ${\mathcal M}:  {\mathbb S}_{h}^r
\rightarrow  {\mathbb S}_{h}^r$ defined via duality by 
\begin{align} 
({\mathcal M}{\mathscr S}_h,\varphi_h)
:&= \frac{\eta}{\tau}({\mathscr S}_h, \varphi_h)  
+ \bigg( \bigg(\frac{i}{\kappa}\nabla+{\bf A}_h^{n+1}\bigg) {\mathscr S}_h \,, 
\bigg(\frac{i}{\kappa}\nabla+{\bf A}_h^{n+1}\bigg) \varphi_h\bigg) \nn \\
&\quad 
+((|{\mathscr S}_h|^{2}-1){\mathscr S}_h, \varphi_h) , 
&&\forall\, \varphi_h \in {\mathbb S}_{h}^r ,
\end{align}
is continuous and monotone, 
i.e.\footnote{The monotonicity makes use of the fact that 
$(|{\mathscr S}_h|^{2}{\mathscr S}_h-|\widetilde{\mathscr S}_h|^{2}\widetilde{\mathscr S}_h,
{\mathscr S}_h-\widetilde{\mathscr S}_h)\geq 0$  
for all ${\mathscr S}_h,\widetilde{\mathscr S}_h \in {\mathbb S}_{h}^r$.}
\begin{align} 
({\mathcal M}{\mathscr S}_h
-{\mathcal M}\widetilde{\mathscr S}_h,{\mathscr S}_h-\widetilde{\mathscr S}_h)
\ge \bigg(\frac{\eta}{\tau}-1\bigg)\|{\mathscr S}_h-\widetilde{\mathscr S}_h\|_{L^2}^2 \, , 
&&\forall\, {\mathscr S}_h,\widetilde{\mathscr S}_h \in {\mathbb S}_{h}^r .
\end{align}
Hence, \cite[Lemma 2.1 and Corollary 2.2 of Chapter 2]{Showalter} 
implies that for any given $f_h\in {\mathbb S}_{h}^r$ 
the equation ${\mathcal M}{\mathscr S}_h=f_h$ 
has a solution ${\mathscr S}_h\in {\mathbb S}_{h}^r$.  
In other words, equation \eqref{FEM1} has a solution $\psi_h^{n+1}\in {\mathbb S}_{h}^r$. The uniqueness of the solution $\psi_h^{n+1}\in {\mathbb S}_{h}^r$ is an obvious consequence of the monotonicity of the operator ${\mathcal M}$. 

Overall, for any given $(\psi_h^{n},{\bf A}_h^{n})\in {\mathbb S}_{h}^{r} \times {\mathbb N}_{h}^{k}$, the system  
\eqref{FEM1}-\eqref{FEM2} has a unique solution 
$(\psi_h^{n+1},\phi_h^{n+1},{\bf A}_h^{n+1})\in {\mathbb S}_{h}^r\times {\mathbb V}_{h}^{k+1} \times {\mathbb N}_{h}^{k}$
when $\tau<\eta$.  
In the rest part of this paper,
we prove the convergence
of the finite element solution.
Some frequently used  
basic lemmas are listed in Section \ref{PreLem}.

\subsection{Preliminary lemmas}\label{PreLem}

The following lemma is concerned with the
approximation properties of the smoothed
projection operators of the finite element spaces \cite{AFW}. 
\begin{lemma}\label{SmoothPr}{\it 
There exist linear projection operators
\begin{align*}
&\widetilde\Pi_h^{\mathbb S}: {\cal L}^1\rightarrow {\mathbb S}_h^r,
\qquad\quad
\widetilde\Pi_h^{\mathbb V}: L^1\rightarrow {\mathbb V}_h^{k+1},
\qquad\quad
\widetilde\Pi_h^{\mathbb N}: {\bf L}^1\rightarrow {\mathbb N}_h^k, 
\end{align*}
which satisfy
\begin{align*}
&\nabla (\widetilde\Pi_h^{\mathbb V}\chi)
=\widetilde\Pi_h^{\mathbb N} \nabla \chi,
&&\forall\,\chi\in W^{1,1}  ,\\
&\|\varphi-\widetilde\Pi_h^{\mathbb S}\varphi\|_{{\cal L}^p}
\leq Ch^{s+3/p-3/q}\|\varphi\|_{{\cal W}^{s,q}}, &&
\forall\,\varphi\in {\cal W}^{s,q},\,\,\,\,\,\,
0\leq s\leq r+1, \\
&\|\chi-\widetilde\Pi_h^{\mathbb V}\chi\|_{L^p}
\leq Ch^{s+3/p-3/q}\|\chi\|_{W^{s,q}}, &&
\forall\,\chi\in W^{s,q},\,\,\,\,\,\,
0\leq s\leq k+2, \\
&\|{\bf a}-\widetilde\Pi_h^{\mathbb N}{\bf a}\|_{{\bf L}^p}
\leq Ch^{s+3/p-3/q}\|{\bf a}\|_{{\bf W}^{s,q}}, 
&&\forall\,{\bf a}\in {\bf W}^{s,q},\,\,\,\,\,\,
0\leq s\leq k+1, 
\end{align*}
for any 
$$
\left\{\begin{array}{ll}
1\leq q\leq p\leq 3/(3/q-s) &\mbox{if}\,\,\,\,\, 0\leq s< 3/q,\\
1\leq q\leq p<\infty &\mbox{if}\,\,\,\,\, s\geq 3/q .
\end{array}
\right. 
$$
}
\end{lemma}

\begin{remark}
The authors of \cite{AFW} (page 66--70) only proved the $L^2$ 
boundedness of the smoothed projection operators.
But their method can also be used to prove
the $L^p$ boundedness without essential change. 
Then Lemma \ref{SmoothPr} is obtained by 
using the Sobolev embedding
$W^{s,q}\hookrightarrow W^{s+3/p-3/q,p}$. 
Although the analysis of \cite{AFW} (page 66--70) only considered polyhedra, 
the extension to curved polyhedra 
is straightforward (as there are no boundary conditions imposed on these 
finite element spaces). 
\end{remark}

It is well known that the solution of the heat equation
\begin{align*} 
\left\{\begin{array}{ll}
\partial_tu-\Delta u =f  & \mbox{in}\,\,\,\Omega ,\\
\nabla u\cdot{\bf n}=0 & \mbox{on}\,\,\,\partial\Omega ,\\
u(x,0)=0, & \mbox{for}\,\,x\in\Omega ,
\end{array}\right. 
\end{align*}
possesses the maximal $L^p$-regularity
(see Corollary 4.d of \cite{Weis1}): 
\begin{align*}
&\|\partial_tu\|_{L^p(0,T;L^q)}
+\|\Delta u\|_{L^p(0,T;L^q)}
\leq C_{p,q}\|f\|_{L^p(0,T;L^q)} ,\qquad 
1<p,q<\infty .
\end{align*}
In this paper, we need to use the maximal $\ell^p$-regularity for 
time-discrete parabolic PDEs, which was proved in \cite[Theorem 3.1]{KLL}. 
\medskip 

\begin{lemma}[Maximal $\ell^p$-regularity]\label{DMPR}{\it  
The solution of the time-discrete PDEs 
\begin{align*}%\label{DMPR2}
\left\{\begin{array}{ll}
D_\tau u^{n+1}-\Delta u^{n+1}=f^{n+1}  & \mbox{in}\,\,\,\Omega ,\\
\nabla u^{n+1}\cdot{\bf n}=0 & \mbox{on}\,\,\,\partial\Omega ,\\
u^0=0,
\end{array}\right. 
\end{align*}
$n=0,1,\cdots$, satisfies
\begin{align*}
&\bigg(\sum_{n=0}^m\tau\|D_\tau u^{n+1}\|_{L^q}^p\bigg)^{\frac{1}{p}}
+\bigg(\sum_{n=0}^m\tau\|\Delta u^{n+1}\|_{L^q}^p\bigg)^{\frac{1}{p}} 
\leq C_{p,q}\bigg(\sum_{n=0}^m\tau\|f^{n+1}\|_{L^q}^p\bigg)^{\frac{1}{p}}
\end{align*}
for any $1<p,q<\infty$ and $m\geq 0$,
where the constant $C_{p,q}$ is independent of 
$\tau$ and $m$. 
} \medskip
\end{lemma}
%\noindent{\it Proof}.$\quad$
%It is well known that the time-continuous PDE 
%\begin{align*}
%\left\{\begin{array}{ll}
%\partial_t u-\Delta u=f & \mbox{in}\,\,\,\Omega ,\\
%\nabla u\cdot{\bf n}=0 & \mbox{on}\,\,\,\partial\Omega ,\\
%u(x,0)=0 & \mbox{for}\,\,\,x\in\Omega , 
%\end{array}\right. 
%\end{align*}
%has the maximal parabolic regularity \cite{Weis}
%\begin{align*}
%\|\partial_tu\|_{L^p(\R_+;L^q)} 
%+\|\Delta u\|_{L^p(\R_+;L^q)} 
%\leq C_{p,q}\|f\|_{L^p(\R_+;L^q)} ,\quad
%\forall\, 1<p,q<\infty .
%\end{align*}
%After a scale transformation $t=s/\tau$ in the time-direction, 
%we see that the problem
%\begin{align*}
%\left\{\begin{array}{ll}
%\partial_s u-\tau \Delta u=f & \mbox{in}\,\,\,\Omega ,\\
%\nabla u\cdot{\bf n}=0 & \mbox{on}\,\,\,\partial\Omega ,\\
%u(x,0)=0 & \mbox{for}\,\,\,x\in\Omega , 
%\end{array}\right. 
%\end{align*}
%also has the maximal parabolic regularity. 
%According to condition (d) of Theorem 1.1 in \cite{Blunck}, 
%the time-discrete problem \eqref{DMPR2} has 
%the discrete maximal parabolic regularity. 
%\endproof

%Some preliminary lemmas which we need to use
%are presented in Section \ref{PLema}. 

We introduce some lemmas in Section \ref{DisHodg} 
on the discrete Hodge 
decomposition, with 
emphasis on the uniform regularity of the discrete
harmonic functions in curved polyhedra. 
A discrete Sobolev embedding inequality 
for functions in the N\'ed\'elec element space
is proved in Section \ref{SobNed}. 
With these mathematical tools, we present estimates
and prove compactness/convergence
of the finite element solution in 
Section \ref{CompFES}.

\subsection{Discrete Hodge decomposition 
and harmonic vector fields}\label{DisHodg}

It is well known that the following Hodge decompositions 
holds (for example, see \cite[decomposition (2.18)]{AFW}) 
\footnote{By identifying the vector fields with the 2-forms, 
in terms of the notation of \cite[decomposition (2.18)]{AFW}, 
we have ${\bf C}(\Omega)\cong {\mathfrak Z}^{*2}$, 
${\bf C}(\Omega)^\perp\cong \mathring{\mathfrak B}^2$, 
${\bf G}(\Omega)\cong {\mathfrak B}^{*2}$ 
and ${\bf X}(\Omega)\cong \mathring{\mathfrak H}^2$.}
\begin{align}\label{Hodge_D_L2}
&{\bf L}^2={\bf C}(\Omega)^\perp\oplus 
{\bf G}(\Omega)\oplus 
{\bf X}(\Omega) , 
\end{align}
where 
\begin{align}
&{\bf C}(\Omega):=\{{\bf u}\in {\bf H}({\rm curl}): \nabla\times{\bf u} =0\}, \\
&{\bf C}(\Omega)^\perp=\{\nabla\times{\bf u}: {\bf u}\in {\bf H}({\rm curl}),\,\,{\bf u}\times{\bf n}=0\} ,\\
&{\bf G}(\Omega):=\{\nabla\omega: \omega\in H^1\} , \\
&{\bf X}(\Omega):=\{{\bf w}\in {\bf H}({\rm curl})\cap{\bf H}({\rm div}):\,
\nabla\times{\bf w}=0, \,\,
\nabla\cdot{\bf w}=0 ,\,\, {\bf w}\cdot{\bf n}=0\,\,\,\mbox{on}
\,\,\,\partial\Omega\}  ,
\end{align}
${\bf C}(\Omega)^\perp$ 
denotes the orthogonal complement of ${\bf C}(\Omega)$ in ${\bf L}^2$, 
and ${\bf X}(\Omega)$ is the space of harmonic vector fields. 

The second type of space of harmonic vector fields 
is defined by \footnote{By identifying the vector fields with the 2-forms, 
in terms of the notation of \cite[definition (2.12)]{AFW}, 
we have $\widetilde {\bf X}(\Omega)={\mathfrak H}^2$.}
\begin{align}
\widetilde {\bf X}(\Omega):=\{{\bf w}\in {\bf H}({\rm curl})\cap{\bf H}({\rm div}): \nabla\times{\bf w}=0,\,\,
\nabla\cdot{\bf w}=0
\,\,\mbox{and}\,\, {\bf w}\times {\bf n}=0
\,\,
\mbox{on}\,\,\partial\Omega\} ,
\end{align}
and we denote 
\begin{align}\label{Def_Space_Y}
\widetilde {\bf Y}(\Omega):=\{{\bf w}\in {\bf H}({\rm curl})\cap{\bf H}({\rm div})\cap 
\widetilde {\bf X}(\Omega)^\perp: {\bf w}\times {\bf n}=0\,\,
\mbox{on}\,\,\partial\Omega \} .
\end{align} 
As a result of \eqref{Hodge_D_L2}, 
any vector field ${\bf v}\in {\bf L}^2$ has the 
Hodge decomposition 
(also see \cite[Appendix]{KY09})  
\begin{align}\label{HodgeDM}
{\bf v} =
\nabla\times {\bf u} +\nabla \omega 
+\sum_{j=1}^{\frak M}\alpha_j{\bf w}_j ,
\end{align}
%\\[-15pt]
%such that\\[-15pt]
%\begin{align}%\label{EL2L3}
%&\|\nabla\times {\bf u}\|_{L^2}
%+\|\nabla \omega\|_{L^2}
%+\sum_{j=1}^{\frak M}|\alpha_j|
%\leq C\|{\bf v}\|_{L^2 } ,
%\end{align}
where ${\bf u}\in \widetilde {\bf Y}(\Omega)$ is the solution of the problem \footnote{By identifying the vector fields with the 2-forms, 
%and identifying the curl operator with the 
%coderivative operator on the 2-forms, 
in terms of the notation of \cite[definition (2.12)]{AFW}, 
we have $\widetilde {\bf X}(\Omega)\cong {\mathfrak H}^2$ 
and $\widetilde {\bf Y}(\Omega)\cong H\Lambda^2(\Omega)\cap 
\mathring H^*\Lambda^2(\Omega)\cap {\mathfrak H}^{2\perp}$. 
Then, by using \cite[Theorem 2.2 on page 23]{AFW} and the Lax--Milgram lemma, 
one can show that the problem \eqref{PDE_hodge_u1}-\eqref{PDE_hodge_u3} has a unique weak solution 
in $\widetilde {\bf Y}(\Omega)$.}
\footnote{\label{FNbd} 
If ${\bf v}\cdot{\bf n}$ is well defined on $\partial\Omega$, then 
the divergence-free part $\nabla\times{\bf u}$ 
satisfies $(\nabla\times{\bf u})\cdot{\bf n}=0$ on $\partial\Omega$, 
due to the boundary conditions in 
\eqref{Eqomega} and \eqref{DHarmF}. 
%This boundary condition is again due to the following 
%fact: ``if a vector field ${\bf u}$ satisfies ${\bf n} \times {\bf u} = 0$ 
%on $\partial\Omega$, then
%$(\nabla\times{\bf u}) \cdot{\bf n}= 0$ on $\partial\Omega$.''
}
\begin{align}
&\nabla\times(\nabla\times {\bf u})=
\nabla\times {\bf v} &&\mbox{in}\,\,\,\Omega,  \label{PDE_hodge_u1}\\
&\nabla\cdot{\bf u}=0 &&\mbox{in}\,\,\,\Omega,\\
&{\bf u}\times {\bf n}=0
&&\mbox{on}\,\,\,\partial\Omega ,  \label{PDE_hodge_u3}
\end{align}
$\omega$ is the solution of 
\begin{align}\label{Eqomega}
&\Delta \omega=\nabla\cdot {\bf v}
&&\mbox{in}\,\,\,\Omega,  \nn\\
&\nabla \omega\cdot{\bf n}={\bf v}\cdot{\bf n}
&&\mbox{on}\,\,\, \partial\Omega,
\end{align}
and ${\bf w}_j=\nabla\varphi_j$, $j=1,2,\cdots,\frak M$, 
form a basis for ${\bf X}(\Omega)$ with $\varphi_j$ being the solution of 
\begin{align}\label{DHarmF}
&\Delta\varphi_j=0 &&\mbox{in}\,\,\,\Omega\backslash\Sigma, \nn\\
&\nabla \varphi_j\cdot{\bf n}=0
&&\mbox{on}\,\,\, \partial\Omega,  \\
&[\nabla \varphi_j\cdot{\bf n}]=0
\quad\mbox{and}\quad [\varphi_j]=\delta_{ij}
&&\mbox{on}\,\,\, \Sigma_i ,\quad i=1,2,\cdots,\frak M,  \nn 
\end{align}
($\delta_{ij}$ denotes the Kronecker symbol). 
The coefficients 
$\alpha_j$, $j=1,\cdots,\frak M$, are given by 
\begin{align}\label{DefAlphaj}
\alpha_j=({\bf v},{\bf w}_j)/\|{\bf w}_j\|_{L^2}^2 . 
\end{align}

\begin{remark}
Although $\varphi_j$ is only defined on 
$\Omega\backslash\Sigma$,
the gradient $\nabla\varphi_j$ has a natural 
extension to be 
a vector field in ${\bf H}({\rm curl,div})$ 
due to the interface conditions.
\end{remark}

To study the regularity of ${\bf w}_j$,
we cite the following
lemma on the regularity of the Poisson
equation in a polyhedral domain. 
This result can be obtained 
by substituting fractional $k$
in Corollary 3.9 of \cite{Dauge} 
(also see page 30 of \cite{Dauge08}
and (23.3) of \cite{Dauge88}). \medskip

\begin{lemma}\label{RegPoiss}
{\it For any given curved polyhedron $\Omega$, there exists a
positive constant $\delta_*>0$ such that the solution
of the Poisson equation
\begin{align*}
\left\{\begin{array}{ll}
-\Delta \varphi=f
&\mbox{in}\,\,\,\Omega, \\
\nabla\varphi\cdot{\bf n}=0
&\mbox{on}\,\,\,\partial\Omega,
\end{array}\right.  
\end{align*}
with the normalization condition 
$\int_\Omega\varphi\d x=0$, satisfies 
\begin{align*} 
\|\varphi\|_{H^{3/2+\alpha}(\Omega)} 
\leq C\|f\|_{H^{-1/2+\alpha}(\Omega)}  
\qquad\mbox{for any $\alpha\in(0,\delta_*]$} .
\end{align*} 
}
\end{lemma}

As a consequence of Lemma \ref{RegPoiss}, 
we have the following result
on the regularity of ${\bf w}_j$, 
(which is also a consequence of Proposition 3.7 of \cite{ABDG}, 
but for self-containedness we include a short proof here). \medskip

\begin{lemma}\label{RegHarV}
{\it 
For any given curved polyhedron $\Omega$, 
there exists a positive constant 
$\delta_*>0$ such that
the harmonic vector fields 
${\bf w}_j$, $j=1,2,\cdots,\frak M$,
are in  ${\bf H}^{1/2+\delta_*}(\Omega)$. \medskip
}
\end{lemma}

\noindent{\it Proof of Lemma \ref{RegHarV}}.$\,\,\,$
Let $\Sigma_{j}'$ 
be a small perturbation of the surfaces
$\Sigma_{j}$ for each $j=1,\cdots,\frak M$, 
such that $\Sigma_j'\cap\Sigma_k=\emptyset$ 
and $\Omega\backslash\Sigma'$ is simply connected
(where $\Sigma'=\cup_{j=1}^{\frak M}\Sigma_j'$). 
Let $D_\Sigma$ and $D_\Sigma'$ be small neighborhoods of
$\Sigma$ and $\Sigma'$, respectively, such that
$\overline D_\Sigma\cap\overline D_\Sigma'=\emptyset$. 

By using Lemma \ref{RegPoiss}   
it is easy to show that the solution
of \eqref{DHarmF} satisfies 
\begin{align*} 
&\varphi_j\in H^{3/2+\delta_*}(\Omega\backslash\overline D_\Sigma) ,
\quad j=1,2,\cdots,\frak M ,
\end{align*}
which implies that
${\bf w}_j=\nabla\varphi_j$, $j=1,\cdots,\frak M$,
are $H^{1/2+\delta_*}$ in the subdomain 
$\Omega\backslash\overline D_\Sigma$. 
Similarly, if we define $\varphi_j'$ as the solution of 
\eqref{DHarmF} with $\Sigma_i$ replaced by $\Sigma_i'$, 
%(and $\Sigma$ replaced by $\Sigma'$),  
%\begin{align*}%\label{DHarmF2}
%&\Delta\varphi_j'=0 &&\mbox{in}\,\,\,\Omega\backslash\Sigma', \nn\\
%&\nabla \varphi_j'\cdot{\bf n}=0
%&&\mbox{on}\,\,\, \partial\Omega,  \\
%&[\nabla \varphi_j'\cdot{\bf n}]=0
%\quad\mbox{and}\quad [\varphi_j']=\delta_{ij}
%&&\mbox{on}\,\,\, \Sigma_i' ,\quad i=1,2,\cdots,\frak M , \nn 
%\end{align*}
then ${\bf w}_j':=\nabla\varphi_j'$, $j=1,\cdots,\frak M$,
also form a basis of ${\bf X}(\Omega)$,
and they are $H^{1/2+\delta_*}$ in the subdomain 
$\Omega\backslash\overline D_\Sigma'$. 
Since ${\bf w}_j$ can be expressed
as linear combinations of ${\bf w}_j'$, it follows that 
${\bf w}_j$ is $H^{1/2+\delta_*}$ in the subdomain 
$\Omega\backslash \overline{D_\Sigma'}\supset\overline D_\Sigma$. 
Therefore, ${\bf w}_j$ is $H^{1/2+\delta_*}$
in the whole domain $\Omega$. 
\endproof \medskip

\begin{definition}{\it 
We define the following 
finite element subspaces of ${\mathbb N}_h^k\subset {\bf H}({\rm curl})$:
\begin{align*}
&{\bf C}_h(\Omega):=\{
{\bf v}_h\in {\mathbb N}_h^k:\,
\nabla\times{\bf v}_h=0 \} , \\
&
{\bf G}_h(\Omega):=\{
\nabla\chi_h:\,
\chi_h\in {\mathbb V}_h^{k+1}\} , \\
&{\bf X}_h(\Omega):=\{
{\bf v}_h\in {\mathbb N}_h^k:\,
\nabla\times {\bf v}_h=0,\,\,
({\bf v}_h,\nabla\chi_h)=0,\,\,
\forall\chi_h\in {\mathbb V}_h^{k+1}\}
\end{align*}
where ${\bf X}_h(\Omega)$
is often referred to as 
the space of discrete harmonic vector fields. \medskip
}
\end{definition}

With the notations above, 
we have the discrete Hodge decomposition 
(page 72 of \cite{AFW}):
\begin{align}\label{DHodgeD}
{\mathbb N}_h^k={\bf C}_h(\Omega)^\perp
\oplus{\bf G}_h(\Omega)\oplus {\bf X}_h(\Omega)  .
\end{align} 
The following lemma  
is concerned with the regularity 
of the discrete harmonic vector
fields.\medskip

\begin{lemma}\label{RegDHarV}
{\it For any given curved polyhedron $\Omega$, 
there exists a positive constant $h_0$ such that
when $h<h_0$ the space ${\bf X}_h(\Omega)$ 
has an orthogonal basis 
$\{{\bf w}_{j,h}$: $j=1,\cdots,\frak M\}$ 
which satisfies  
\begin{align}%\label{EL2L30}
&\sum_{j=1}^{\frak M}\|{\bf w}_{j,h}\|_{L^{3+\delta}}\leq C 
\qquad\mbox{and}\qquad
\sum_{j=1}^{\frak M}\|{\bf w}_{j,h}-{\bf w}_j\|_{L^{3+\delta}}
\rightarrow 0 \quad\mbox{as}\,\,\, h\rightarrow 0  , 
\end{align}
for any $0<\delta<3\delta_*/(1-\delta_*)$, where  
$\delta_*$ is given by Lemma \ref{RegHarV}.
}
\medskip
\end{lemma}

\noindent{\it Proof of Lemma \ref{RegDHarV}}.$\,\,\,$
If ${\bf v}_h\in {\bf X}_h(\Omega)$,
then $\nabla\times {\bf v}_h=0$ and so 
the Hodge decomposition \eqref{HodgeDM} implies 
\begin{align*} {\bf v}_h =
\nabla \omega 
+\sum_{j=1}^{\frak M}\alpha_j{\bf w}_j .
\end{align*}
Using the commuting property 
of the smoothed projection operator (Lemma \ref{SmoothPr}) 
we derive 
\begin{align} \label{DPvh}
{\bf v}_h 
&= \widetilde\Pi_h^{\mathbb N}\nabla \omega 
+\sum_{j=1}^{\frak M}\alpha_j\widetilde\Pi_h^{\mathbb N} {\bf w}_j 
= \nabla \widetilde\Pi_h^{\mathbb V}\omega 
+\sum_{j=1}^{\frak M}\alpha_j\widetilde\Pi_h^{\mathbb N} {\bf w}_j 
=: \nabla \omega_h 
+\sum_{j=1}^{\frak M}\alpha_j\widetilde\Pi_h^{\mathbb N} {\bf w}_j ,
\end{align}
where we have defined
$\omega_h:=\widetilde\Pi_h^{\mathbb V}\omega$
to simplify the notation.
Since any ${\bf v}_h\in {\bf X}_h(\Omega)$ 
satisfies $({\bf v}_h,\nabla\chi_h)=0$ for all 
$\chi_h\in {\mathbb V}_h^{k+1}$,
it follows that 
\begin{align*} 
(\nabla \omega_h,\nabla\chi_h)
=-\sum_{j=1}^{\frak M}\alpha_j
(\widetilde\Pi_h^{\mathbb N} {\bf w}_j,\nabla\chi_h)
=\sum_{j=1}^{\frak M}\alpha_j
({\bf w}_j-\widetilde\Pi_h^{\mathbb N} {\bf w}_j,\nabla\chi_h),
\quad\forall\, \chi_h\in {\mathbb V}_h^{k+1} .
\end{align*}

If we define $\omega_{j,h}\in{\mathbb V}_h^{k+1}$ 
(with the normalization $\int_\Omega\omega_{j,h}\d x=0$) 
as the finite element solution of 
\begin{align}\label{omegajh}
(\nabla \omega_{j,h},\nabla\chi_h)
=({\bf w}_j-\widetilde\Pi_h^{\mathbb N} {\bf w}_j,\nabla\chi_h),
\quad\forall\, \chi_h\in {\mathbb V}_h^{k+1} ,
\end{align}
then we have
$
\omega_{h}=\sum_{j=1}^{\frak M}
\alpha_j \omega_{j,h}  
+{\rm const} .
$
Substituting this into \eqref{DPvh}, 
we obtain
\begin{align*} 
{\bf v}_h
=\sum_{j=1}^{\frak M}\alpha_j(\nabla \omega_{j,h} 
+ \widetilde\Pi_h^{\mathbb N} {\bf w}_j) .
\end{align*}
We see that any vector field in ${\bf X}_h(\Omega)$
can be expressed as a linear combination
of  
\begin{align} 
{\bf w}_{j,h}:=\nabla\omega_{j,h} 
+\widetilde\Pi_h^{\mathbb N} {\bf w}_j ,\quad
j=1,\cdots,\frak M .
\end{align} 
The vector fields ${\bf w}_{j,h}$, $j=1,\cdots,\frak M$, 
must form a basis for ${\bf X}_h(\Omega)$
if they are linearly independent. 
Indeed, by substituting $\chi_h=\omega_{j,h}$
into \eqref{omegajh}, we obtain
\begin{align*} 
\|\nabla \omega_{j,h}\|_{L^2}\leq
C\|{\bf w}_j-\widetilde\Pi_h^{\mathbb N} {\bf w}_j\|_{L^2}
\leq Ch^{1/2+\delta_*}\|{\bf w}_j\|_{H^{1/2+\delta_*}} .
%\leq Ch^{1/2+\delta_*} .
\end{align*} 
Using the inverse inequality, we 
see that for $\delta<3\delta_*/(1-\delta_*)$ there holds 
\begin{align*} 
\|\nabla \omega_{j,h}\|_{L^{3+\delta}}\leq
Ch^{-1/2-\delta/(3+\delta)}
\|\nabla \omega_{j,h}\|_{L^2}
\leq Ch^{\delta_*-\delta/(3+\delta)} \rightarrow 0
\quad\mbox{as}\,\,\, h\rightarrow 0 .
\end{align*} 
%Since
%\begin{align*} 
%\|\widetilde\Pi_h^{\mathbb N} {\bf w}_j-{\bf w}_j\|_{L^{3+\delta}}
%\leq Ch^{\delta_*-\delta/(3+\delta)} \|{\bf w}_j\|_{H^{1/2+\delta_*}} 
%\leq Ch^{\delta_*-\delta/(3+\delta)} 
%\rightarrow 0
%\quad\mbox{as}\,\,\, h\rightarrow 0 ,
%\end{align*} 
Using Lemma \ref{SmoothPr} and Lemma \ref{RegHarV}, we have 
\begin{align*} 
\|{\bf w}_{j,h}-{\bf w}_j\|_{L^{3+\delta}}
&\leq
\|\nabla \omega_{j,h}\|_{L^{3+\delta}}
+\|\widetilde\Pi_h^{\mathbb N} {\bf w}_j-{\bf w}_j\|_{L^{3+\delta}}\\
&\leq %Ch^{\delta_*-\delta/(3+\delta)} +
Ch^{\delta_*-\delta/(3+\delta)}\|{\bf w}_j\|_{H^{1/2+\delta_*}} 
 \rightarrow 0\quad \mbox{as}\,\,\,h\rightarrow 0,
&& j=1,\cdots,\frak M. 
\end{align*} 
Since ${\bf w}_j$, $j=1,\cdots,\frak M$,
are linearly independent and ${\bf w}_{j,h}$ converges to
${\bf w}_{j}$, 
there exists a positive constant $h_0$ such that
${\bf w}_{j,h}$, $j=1,\cdots,\frak M$,
are also linearly independent when $h<h_0$. 

A Gram--Schmidt orthogonalization
process gives an orthogonal
basis which still converges to 
the basis of ${\bf X}(\Omega)$
in ${\bf L}^{3+\delta}$. 
The proof of Lemma \ref{RegDHarV} is complete.
\endproof

\subsection{A discrete Sobolev embedding inequality 
for the N\'ed\'elec element space}\label{SobNed}
\indent
\begin{definition}\label{DisDiv}{\it 
For any given ${\bf a}_h\in {\mathbb N}_{h}^{k}$, 
the unique function $\zeta_h\in {\mathbb V}_h^{k+1}$ satisfying 
$$ 
(\zeta_h, \chi_h)=-({\bf a}_h,\nabla\chi_h) ,
\quad\forall\, \chi_h\in {\mathbb V}_h^{k+1} ,
$$ 
is called the  
discrete divergence of ${\bf a}_h$,  
denoted by 
$\zeta_h:=\nabla_h^{\mathbb N}\cdot{\bf a}_h$. 
The discrete analogue of the ${\bf H}({\rm curl,div})$ norm 
is defined as
\begin{align}
\|{\bf a}_h\|_{{\bf H}_h({\rm curl,div})}
:=\|{\bf a}_h\|_{L^2}+\|\nabla\times{\bf a}_h\|_{L^2} 
+\|\nabla_h^{\mathbb N}\cdot{\bf a}_h\|_{L^2}  .
\end{align}
}
\end{definition}

\begin{lemma}\label{SobolevD}
{\it For any given curved polyhedron $\Omega$, 
there exist positive constants $h_0$, $\delta$
and $C$ such that if the set of functions 
$\{{\bf a}_h\in {\mathbb N}_{h}^{k}:\, h>0\}$ is bounded in
the norm $\|\cdot\|_{{\bf H}_h({\rm curl,div})}$,
then it is compact 
in ${\bf L}^{3+\delta}$, 
and 
\begin{align}\label{aaa1}
\|{\bf a}_h\|_{L^{3+\delta}}
\leq C\|{\bf a}_h\|_{{\bf H}_h({\rm curl,div})} 
\quad\mbox{when}\,\,\, h<h_0.
\end{align}
}
\end{lemma}

\noindent{\it Proof of Lemma \ref{SobolevD}}.$\,\,\,$ 
The discrete Hodge 
decomposition \eqref{DHodgeD} implies 
\begin{align}
{\bf a}_h =
{\bf c}_h + \nabla\theta_h 
+ \sum_{j=1}^{\frak M}\alpha_{j,h}{\bf w}_{j,h} ,
\end{align}
where ${\bf c}_h\in {\bf C}_h(\Omega)^\perp$, 
$\theta_h\in {\mathbb V}_h^{k+1}$
and ${\bf w}_{j,h}$, $j=1,\cdots,\frak M$,
are the basis functions of ${\bf X}_h(\Omega)$
given in Lemma \ref{RegDHarV}. 
We shall prove that the three functions 
are all compact in ${\bf L}^{3+\delta}(\Omega)$. 
%\begin{align}
%\|{\bf c}_h\|_{L^{3+\delta} }
%\leq C\|\nabla\times {\bf a}_h\|_{L^2 }  
%\end{align}
%and
%\begin{align}
%\|\nabla\theta_h\|_{L^{3+\delta} }
%\leq C\|\zeta_h\|_{L^2 } 
%\end{align}
%separately. 

Firstly, consider the continuous Hodge decomposition of 
${\bf a}_h$ (see \eqref{HodgeDM})  
\begin{align}
{\bf a}_h =
\nabla\times {\bf u}^h +\nabla \omega^h 
+\sum_{j=1}^{\frak M}\alpha_j^h{\bf w}_j ,
\end{align}
where ${\bf u}^h\in \widetilde {\bf Y}(\Omega)$ 
is the solution of the PDE problem 
\footnote{See \eqref{Def_Space_Y}-\eqref{PDE_hodge_u3} for the definition of the space $\widetilde {\bf Y}(\Omega)$.}
\begin{align*} 
&\nabla\times(\nabla\times {\bf u}^h)=\nabla\times{\bf a}_h ,
&&\mbox{in}\,\,\,\Omega, \nn\\
&\nabla\cdot{\bf u}^h=0 , &&\mbox{in}\,\,\,\Omega,  \\
&{\bf u}^h\times{\bf n}=0
&&\mbox{on}\,\,\,\partial\Omega .  \nn
\end{align*} 
Hence, the vector field 
${\bf c}^h:=\nabla\times {\bf u}^h\in {\bf C}(\Omega)^\perp$ 
is the divergence-free part 
of ${\bf a}_h$,  which satisfies 
${\bf c}^h\cdot{\bf n}=0$ \footnote{See\, footnote\, \ref{FNbd} 
on this boundary condition.}
and the basic energy inequality
\begin{align}\label{ccc1}
\|{\bf c}^h\|_{{\bf H}({\rm curl,div})}
\leq C\|\nabla\times{\bf a}_h\|_{L^2 } .
%\quad \mbox{for some $s\in(1/2,1]$}.
\end{align} 
Since ${\bf H}({\rm curl,div})\hookrightarrow
{\bf H}^{1/2+\delta_*}(\Omega)$ for some $\delta_*>0$ 
\footnote{This is a immediate consequence of 
Lemma \ref{RegPoiss} and the following 
decomposition proved in \cite{BS87}: 
$$
{\bf H}({\rm curl,div})=
{\bf H}^1+ \{\nabla\varphi:\varphi\in H^1,\,\,
\Delta \varphi\in L^2,\,\,\nabla\varphi\cdot{\bf n}
=0\,\,\mbox{on}\,\,\partial\Omega\}. 
$$ }
and ${\bf H}^{1/2+\delta_*}(\Omega)$ is 
compactly embeddded into ${\bf L}^{3+\delta}(\Omega)$
for $\delta<3\delta_*/(1-\delta_*)$,
it follows that the set 
$\{{\bf c}^h:\, h>0\}$ is 
compact in ${\bf L}^{3+\delta}(\Omega)$. 

Since 
$$
\nabla\times({\bf c}^h-{\bf c}_h)
=\nabla\times {\bf c}^h-\nabla\times {\bf c}_h 
=\nabla\times {\bf a}_h
- \nabla\times{\bf a}_h =0 ,
$$ 
it follows from \cite[Theorem 5.11 on page 74]{AFW} that 
\footnote{By identifying the vector fields with the 
1-forms, in terms of the notation of \cite[Theorem 5.11 on page 74]{AFW}, 
we have ${\bf C}(\Omega)\cong {\mathfrak Z}^{1}$ 
and ${\bf C}(\Omega)^\perp\cong {\mathfrak Z}^{1\perp}$.} 
\begin{align*} 
\|\widetilde\Pi_h^{\mathbb N}{\bf c}^h-{\bf c}_h\|_{L^2}
\leq C\|{\bf c}^h\|_{H^{1/2+\delta_*}}h^{1/2+\delta_*}
\leq C\|\nabla\times{\bf a}_h\|_{L^2} h^{1/2+\delta_*} ,
\end{align*} 
and by using the inverse inequality we further derive 
\begin{align*} 
\|\widetilde\Pi_h^{\mathbb N}{\bf c}^h-{\bf c}_h\|_{L^{3+\delta}}
&\leq Ch^{-1/2-\delta/(3+\delta)}\|\widetilde\Pi_h^{\mathbb N}{\bf c}^h
-{\bf c}_h\|_{L^2} 
\leq C\|\nabla\times{\bf a}_h\|_{L^2}
h^{\delta_*-\delta/(3+\delta)} .
\end{align*}
Since $\delta_*-\delta/(3+\delta)>0$
when $\delta<3\delta_*/(1-\delta_*)$, 
by using Lemma \ref{SmoothPr} we have   
\begin{align} \label{ccc2}
\|{\bf c}^h-{\bf c}_h\|_{L^{3+\delta} }
&\leq \|\widetilde\Pi_h^{\mathbb N}{\bf c}^h-{\bf c}^h\|_{L^{3+\delta} }
+\|\widetilde\Pi_h^{\mathbb N}{\bf c}^h-{\bf c}_h\|_{L^{3+\delta} }\nn\\
&\leq 
C\|{\bf c}^h\|_{H^{1/2+\delta_*}(\Omega)}h^{\delta_*-\delta/(3+\delta)}
+C\|\nabla\times{\bf a}_h\|_{L^2 }
h^{\delta_*-\delta/(3+\delta)}\nn\\
&
\leq C\|\nabla\times{\bf a}_h\|_{L^2 }
h^{\delta_*-\delta/(3+\delta)} \rightarrow 0
\qquad\mbox{as}\,\,\, h\rightarrow 0 .
\end{align}
Since $\{{\bf c}^h:\, h>0\}$ 
is compact in ${\bf L}^{3+\delta}(\Omega)$
and $\|{\bf c}^h-{\bf c}_h\|_{L^{3+\delta} }\rightarrow 0$
as $h\rightarrow 0$, it follows that $\{{\bf c}_h:\, h>0\}$ 
is also compact in ${\bf L}^{3+\delta}(\Omega)$. 

Secondly, we let
$\zeta_h=\nabla_h^{\mathbb N}\cdot{\bf a}_h$
in the sense of Definition \ref{DisDiv}.
Due to the orthogonality of 
${\bf c}_h$ and ${\bf w}_{j,h}$ with $\nabla\chi_h$, we have 
\begin{align*}
(\nabla\theta_h,\nabla\chi_h)
=({\bf a}_h,\nabla \chi_h) 
=-(\zeta_h, \chi_h) , \quad\forall\,\chi_h\in {\mathbb V}_h^{k+1} .
\end{align*}
Let $\theta^h$ be the solution of the PDE problem
\begin{align*}
&\Delta\theta^h=\zeta_h &&\mbox{in}\,\,\,\Omega, \nn\\
&\nabla\theta^h\cdot{\bf n}=0 &&\mbox{on}\,\,\,\partial \Omega  ,
\end{align*}
which satisfies (using Lemma \ref{RegPoiss}) 
\begin{align}\label{ccc3}
\|\theta^h\|_{H^{3/2+\delta_*}(\Omega)} 
\leq C\|\zeta_h\|_{L^2 } 
\quad\mbox{for some $\delta_*>0$} .
\end{align}
Hence, the set $\{\nabla\theta^h:\, h>0\}$
is bounded in ${\bf H}^{1/2+\delta_*}(\Omega)$, which is compactly embedded
into ${\bf L}^{3+\delta}(\Omega)$
for $\delta< 3\delta_*/(1-\delta_*)$. 
Moreover, according to the definition of $\theta^h$, 
we have 
\begin{align*}
(\nabla(\theta^h-\theta_h),\nabla\chi_h)
=0, \quad\forall\,\chi_h\in {\mathbb V}_h^{k+1}  .
\end{align*}
By substituting $\chi_h=\widetilde\Pi_h^{\mathbb V}\theta^h-\theta_h$
into the last equation, we obtain 
\begin{align*}
\|\nabla(\widetilde\Pi_h^{\mathbb V}\theta^h-\theta_h)\|_{L^2 }
\leq C\|\theta^h\|_{H^{3/2+\delta_*}(\Omega)}h^{1/2+\delta_*}
\leq C\|\zeta_h\|_{L^2 }h^{1/2+\delta_*} .
\end{align*}
Again, by using the inverse inequality we derive
\begin{align*}
\|\nabla(\widetilde\Pi_h^{\mathbb V}\theta^h-\theta_h)\|_{L^{3+\delta} }
&\leq Ch^{-1/2-\delta/(3+\delta)}
\|\nabla(\widetilde\Pi_h^{\mathbb V}\theta^h-\theta_h)
\|_{L^2 } 
\leq C\|\zeta_h\|_{L^2 }h^{\delta_*-\delta/(3+\delta)} . 
\end{align*} 
In view of Lemma \ref{SmoothPr}, we have 
\begin{align} \label{ccc4}
\|\nabla\theta^h-\nabla\theta_h\|_{L^{3+\delta} }
&\leq \|\nabla(\theta^h-\widetilde\Pi_h^{\mathbb V}\theta^h)\|_{L^{3+\delta} }
+\|\nabla(\widetilde\Pi_h^{\mathbb V}\theta^h-\theta_h)\|_{L^{3+\delta} } \nn\\
&\leq C\|\theta^h\|_{H^{3/2+\delta_*}(\Omega)}
h^{\delta_*-\delta/(3+\delta)}
+C\|\zeta_h\|_{L^2 }h^{\delta_*-\delta/(3+\delta)}\nn\\
&
\leq C\|\zeta_h\|_{L^2 }h^{\delta_*-\delta/(3+\delta)} 
\rightarrow 0\quad\mbox{as}\,\,\, h\rightarrow 0 .
\end{align}
Therefore, the set of functions $\{\nabla\theta_h:\, h>0\}$ 
is compact in ${\bf L}^{3+\delta}(\Omega)$. 

Finally, we note that
\begin{align}\label{Alphaj2}
|\alpha_{j,h}|=|({\bf a}_h,{\bf w}_{j,h})|/\|{\bf w}_{j,h}\|_{L^2}^2
\leq C\|{\bf a}_h\|_{L^2} 
\leq C\|{\bf a}_h\|_{{\bf H}_h({\rm curl,div})} ,
\quad j=1,\cdots,\frak M. 
\end{align}
Therefore, the set of numbers 
$\{\alpha_{j,h}:\, h>0\}$, are compact. 
Since ${\bf w}_{j,h}$ converges to ${\bf w}_j$
in ${\bf L}^{3+\delta}(\Omega)$ 
(see Lemma \ref{RegDHarV}),
it follows that 
$\big\{\sum_{j=1}^{\frak M}\alpha_{j,h}{\bf w}_{j,h}:\, h>0\big\}$
is compact in ${\bf L}^{3+\delta}(\Omega)$. 

Overall, we have proved that
${\bf c}_h$, $\nabla\theta_h$ and 
$\sum_{j=1}^{\frak M}\alpha_{j,h}{\bf w}_{j,h}$
are all compact in ${\bf L}^{3+\delta}(\Omega)$. 
The inequalities \eqref{ccc1} and \eqref{ccc4}-\eqref{Alphaj2} imply \eqref{aaa1}.  
The proof of Lemma \ref{SobolevD} is complete. \,\endproof\medskip

\begin{remark}
If the domain $\Omega$ is smooth or convex,
then a similar proof yields 
\begin{align}\label{aaa-2}
\|{\bf a}_h\|_{L^{6}}
\leq C\|{\bf a}_h\|_{{\bf H}_h({\rm curl,div})} .
\end{align}
\end{remark}

\subsection{Uniform estimates
of the finite element solution}\label{SecUnif}

In this subsection we prove the following
lemma. 
\begin{lemma}\label{UniFEst}
{\it There exist positive constants $\tau_0\in(0,\eta/2)$,
$q>3$ and $C$  
such that when $\tau<\tau_0$
the finite element solution satisfies 
\begin{align}\label{LUniFEst}
&\max_{0\leq n\leq N-1}
\big(\|\psi^{n+1}_h\|_{H^1}+\|{\bf A}^{n+1}_h\|_{L^q} 
+\|\phi_h^{n+1}\|_{L^2}+
\|\nabla\times{\bf A}_h^{n+1}\|_{L^2}\big)  \nn \\
&+\sum_{n=0}^{N-1}\tau \big(\|D_\tau\psi_h^{n+1}\|_{L^2}^2 
+\|D_\tau{\bf A}_h^{n+1}\|_{L^2}^2 \big) \nn\\
&+\sum_{n=0}^{N-1}\tau \big(\|\psi_{h}^{n+1}\|_{W^{1,q}}^2
+\|\phi_{h}^{n+1}\|_{H^1}^2
+\|D_\tau\phi_{h}^{n+1}\|_{(H^1)'}^2 \big)
\leq C .  
\end{align}
}
\end{lemma}

\noindent{\it Proof of Lemma \ref{UniFEst}}.$\quad$
We shall prove the following inequality by 
mathematical induction:
\begin{align}\label{MathInd}
\|\psi_h^n\|_{L^\infty}\leq \tau^{-1/2} .
\end{align}
Since $|\psi_h^0|\leq 1$,
it follows that \eqref{MathInd} holds for $n=0$
when $\tau<1$.
In the following, we assume that the inequality holds for 
$0\leq n\leq m\leq N-1$ and prove that
it also holds for $n=m+1$. 
The generic constant $C$ of this subsection will
be independent of $h$, $\tau$ and $m$.

Under the induction assumption above, 
from \eqref{FEMEGDecay0} we see that 
%(using the discrete Gronwall's inequality) 
\begin{align*}%\label{FEEnerg}
&\max_{0\leq n\leq m}{\cal G}_h^{n+1} 
+\sum_{n=0}^m\tau \int_\Omega 
\bigg(\frac{\eta-\tau}{2} |D_\tau\psi_h^{n+1}|^2
+\frac{1}{2}|D_\tau{\bf A}_h^{n+1}|^2 \bigg) \d x
\leq C ,
%\leq e^{\eta\kappa^2T}{\cal G}^{**}(\psi_h^0,{\bf A}_h^0,{\bf B}_h^0)
\end{align*}
which implies 
\begin{align}\label{FreeEngUni}
&\max_{0\leq n\leq m}
\bigg(\bigg\|\frac{i}{\kappa}\nabla\psi^{n+1}_h 
+ \mathbf{A}^{n+1}_h\psi^{n+1}_h\bigg\|_{L^2}
+\|\psi_h^{n+1}\|_{L^4} \bigg) \nn\\
&
+\max_{0\leq n\leq m}
\bigg(\|\phi_h^{n+1}\|_{L^2}+
\|\nabla\times{\bf A}_h^{n+1}\|_{L^2}
+\|{\bf A}_h^{n+1}\|_{L^2}\bigg)\nn\\
&+\sum_{n=0}^m\tau \big(\|D_\tau\psi_h^{n+1}\|_{L^2}^2 
+\|D_\tau{\bf A}_h^{n+1}\|_{L^2}^2\big)
\leq C . 
%\leq e^{\eta\kappa^2T}{\cal G}^{**}(\psi_h^0,{\bf A}_h^0,{\bf B}_h^0)
\end{align}
We assume $0\leq n\leq m$ below
if there is no explicit mention of the range of $n$,
and let $\ell^p_{m}(W^{l,q})$ denote the space 
of sequences $(v_n)_{n=0}^m$, with 
$v_n\in W^{l,q}$, 
equipped with the following norm: 
$$
\|(v_{n})_{n=0}^m\|_{\ell^p(W^{l,q})}
:=\left\{\begin{array}{ll}
\displaystyle\bigg( \sum_{n=0}^m\tau 
\|v_{n}\|_{W^{l,q}}^p\bigg)^{\frac{1}{p}} 
&\quad\mbox{if}\,\,\, 1\leq p<\infty 
\,\,\,\mbox{and}\,\,\, 1\leq q\leq \infty ,\\[15pt]
\displaystyle
\max_{0\leq n\leq m}\|v_{n}\|_{W^{l,q}}
&\quad\mbox{if}\,\,\,  p=\infty \,\,\,\mbox{and}\,\,\, 1\leq q\leq \infty .
\end{array}
\right. 
$$

In view of \eqref{FEM3}, Lemma \ref{SobolevD} 
implies the existence of $q>3$ 
such that 
\begin{align}\label{AhLinftLq} 
\max_{0\leq n\leq m}\|{\bf A}^{n+1}_h\|_{L^q}
&\leq C\max_{0\leq n\leq m} (\|\phi_h^{n+1}\|_{L^2}
+\|\nabla\times{\bf A}^{n+1}_h\|_{L^2} 
+\|{\bf A}_h^{n+1}\|_{L^2} ) 
 \leq C .
\end{align}
Let $\bar q<6$ be the number satisfying 
$1/q+1/\bar q=1/2$. 
By using H\"older's inequality we derive 
\begin{align*} 
\|{\bf A}^{n+1}_h\psi_h^{n+1}\|_{L^2} 
&\leq C \|{\bf A}^{n+1}_h\|_{L^q}\|\psi_h^{n+1}\|_{L^{\bar q}} 
\leq C\|\psi_h^{n+1}\|_{L^{\bar q}} 
\leq \epsilon\|\nabla\psi_h^{n+1}\|_{L^2} 
+C_\epsilon\|\psi_h^{n+1}\|_{L^2} ,
\end{align*}
where we have also used the 
interpolation inequality
\begin{align*} 
\|\psi_h^{n+1}\|_{L^{\bar q}} 
%&\leq \|\psi_h^{n+1}\|_{L^2}^{3/\bar q-1/2}
%\|\psi_h^{n+1}\|_{L^6}^{3/2-3/\bar q} \\
&\leq C\|\psi_h^{n+1}\|_{L^2}^{3/\bar q-1/2}
\|\psi_h^{n+1}\|_{H^1}^{3/2-3/\bar q}
\leq \epsilon\|\nabla\psi_h^{n+1}\|_{L^2} 
+C_\epsilon\|\psi_h^{n+1}\|_{L^2} ,
\quad
\forall\, \epsilon\in(0,1) .
\end{align*} 
As a consequence, we have
\begin{align*}
\|\nabla\psi_h^{n+1}\|_{L^2}
&\leq \bigg\|\frac{i}{\kappa}\nabla\psi^{n+1}_h 
+ \mathbf{A}^{n+1}_h\psi^{n+1}_h\bigg\|_{L^2}
+\|{\bf A}^{n+1}_h\psi_h^{n+1}\|_{L^2} \nn \\
&\leq 
\bigg\|\frac{i}{\kappa}\nabla\psi^{n+1}_h 
+ \mathbf{A}^{n+1}_h\psi^{n+1}_h\bigg\|_{L^2}
+ \epsilon\|\nabla\psi_h^{n+1}\|_{L^2} 
+C_\epsilon\|\psi_h^{n+1}\|_{L^2} ,
\end{align*}
which further reduces to
(by choosing $\epsilon=1/2$)
\begin{align}\label{psihLinfH1}
\max_{0\leq n\leq m}\|\nabla\psi_h^{n+1}\|_{L^2}
&\leq 
C\max_{0\leq n\leq m}\bigg\|\frac{i}{\kappa}\nabla\psi^{n+1}_h 
+ \mathbf{A}^{n+1}_h\psi^{n+1}_h\bigg\|_{L^2}
+C\max_{0\leq n\leq m}\|\psi_h^{n+1}\|_{L^2}^2 \leq C .
\end{align}

To estimate $\|\psi_h^{n+1}\|_{L^\infty}$, 
we need the
following lemma. \medskip

\begin{lemma}\label{LemmW1qD}
{\it There exists a positive constant $q_0\in(3,4]$ such that 
for $3<q<q_0$ the finite element solution 
$\psi_h^{n+1}\in {\mathbb S}_h^r$, 
$n=0,1,\cdots,m$, 
of the equation
\begin{align}\label{EQDPSI}
\eta D_\tau\psi_h^{n+1}
-\frac{1}{\kappa^2}\Delta_h\psi_h^{n+1}
=f_h^{n+1}
\end{align}
satisfies 
\begin{align}\label{LemW1qD}
\|(\psi_h^{n+1})_{n=0}^m\|_{\ell^2(W^{1,q+\delta_q})}\leq 
C\|(f_h^{n+1})_{n=0}^m\|_{\ell^2(L^{q/2})}
+C\|\psi_h^0\|_{H^1}
\quad\mbox{
for some $\delta_q>0$. }
\end{align} 
}
\end{lemma}
\noindent{\it Proof of Lemma \ref{LemmW1qD}}.$\quad$
Let $\theta^{n+1}$ be the solution of the PDE
problem 
\begin{align}\label{EQDPSI2}
\left\{\begin{array}{ll}
\displaystyle
\eta D_\tau\theta^{n+1}-\frac{1}{\kappa^2}\Delta\theta^{n+1}
=f_h^{n+1} &\mbox{in}\,\,\,\Omega,\\[10pt]
\nabla\theta^{n+1}\cdot{\bf n}=0 
&\mbox{on}\,\,\,\partial\Omega ,\\[5pt]
\theta^0=\psi_h^0 .
\end{array}\right. 
\end{align}
The function $\theta^{n+1}$ can further 
be decomposed as
$\theta^{n+1}=\widehat\theta^{n+1}
+\widetilde\theta^{n+1}$, which are solutions of
\begin{align*}
\left\{\begin{array}{ll}
\displaystyle
\eta D_\tau\widehat\theta^{n+1}-\frac{1}{\kappa^2}\Delta\widehat\theta^{n+1}
=f_h^{n+1} &\mbox{in}\,\,\,\Omega,\\[10pt]
\nabla\widehat\theta^{n+1}\cdot{\bf n}=0 
&\mbox{on}\,\,\,\partial\Omega ,\\[5pt]
\widehat\theta^0=0 .
\end{array}\right. 
\quad\mbox{and}\quad 
\left\{\begin{array}{ll}
\displaystyle
\eta D_\tau\widetilde\theta^{n+1}-\frac{1}{\kappa^2}\Delta\widetilde\theta^{n+1}
=0 &\mbox{in}\,\,\,\Omega,\\[10pt]
\nabla\widetilde\theta^{n+1}\cdot{\bf n}=0 
&\mbox{on}\,\,\,\partial\Omega ,\\[5pt]
\widetilde\theta^0=\psi_h^0 ,
\end{array}\right. 
\end{align*}
respectively. 
The solution $\widehat\theta^{n+1}$
satisfies (see Lemma \ref{DMPR}) 
\begin{align*}
\|(D_\tau\widehat\theta^{n+1})_{n=0}^m\|_{\ell^2(L^{q/2})}
+\|(\Delta\widehat\theta^{n+1})_{n=0}^m\|_{\ell^2(L^{q/2})} \leq 
C\|(f_h^{n+1})_{n=0}^m\|_{\ell^2(L^{q/2})} ,
\quad \forall\, 2<q<\infty, 
\end{align*}
and $\widetilde\theta^{n+1}$ satisfies the standard energy
estimate 
\begin{align*}
\|(D_\tau\widetilde\theta^{n+1})_{n=0}^m\|_{\ell^2(L^2)}
+\|(\Delta\widetilde\theta^{n+1})_{n=0}^m\|_{\ell^2(L^2)} \leq 
C\|\widetilde\theta^0\|_{H^1} .
\end{align*}
In view of the last two inequalities, for any $2<q\leq 4$
we have
\begin{align}\label{thL2Lq2}
\|(D_\tau\theta^{n+1})_{n=0}^m\|_{\ell^2(L^{q/2})}
+\|(\Delta\theta^{n+1})_{n=0}^m\|_{\ell^2(L^{q/2})} \leq 
C\|(f_h^{n+1})_{n=0}^m\|_{\ell^2(L^{q/2})} 
+C\|\psi_h^0\|_{H^1} .
\end{align} 
If we define $\overline\theta^{n+1}:=
\frac{1}{|\Omega|}\int_\Omega\theta^{n+1}\d x$
as the average of $\theta^{n+1}$ over $\Omega$, then 
Lemma \ref{RegPoiss} implies 
\begin{align*} 
\|(\theta^{n+1}-\overline\theta^{n+1})_{n=0}^m \|_{\ell^2(H^{3/2+\alpha})}  
\leq  
C\|(\Delta\theta^{n+1})_{n=0}^m\|_{\ell^2(H^{-1/2+\alpha}) } 
\end{align*} 
for any $0<\alpha<\min(\delta_*,\frac{1}{2})$. 
The last inequality implies 
\begin{align}\label{thL2Lq3}
\|(\theta^{n+1})_{n=0}^m\|_{\ell^2(H^{3/2+\alpha})}  
\leq  
C\|(\Delta\theta^{n+1})_{n=0}^m\|_{\ell^2(H^{-1/2+\alpha}) } 
+C\|(\theta^{n+1})_{n=0}^m\|_{\ell^2(L^1)}  .
\end{align} 
For any 
\begin{align}\label{Defalpha}
3<q=6/(2-\alpha)<\min(6/(2-\delta_*),4)   ,
\end{align} 
the Sobolev embedding 
$L^{q/2}\hookrightarrow H^{-1/2+\alpha}$ 
and \eqref{thL2Lq2}-\eqref{thL2Lq3} imply 
\begin{align*} 
\|(\theta^{n+1})_{n=0}^m\|_{\ell^2(H^{3/2+\alpha})}  
&\leq C\|(\Delta\theta^{n+1})_{n=0}^m\|_{\ell^2(H^{-1/2+\alpha}) }
+C\|(\theta^{n+1})_{n=0}^m\|_{\ell^2(L^1)}  \nn\\
&\leq 
C\|(\Delta\theta^{n+1})_{n=0}^m\|_{\ell^2(L^{q/2})}
+C\|(D_\tau\theta^{n+1})_{n=0}^m\|_{\ell^2(L^1)}
+C\|\theta^0\|_{L^1} \nn\\
&\leq 
C\|(\Delta\theta^{n+1})_{n=0}^m\|_{\ell^2(L^{q/2})}
+C\|(D_\tau\theta^{n+1})_{n=0}^m\|_{\ell^2(L^{q/2})}
+C\|\theta^0\|_{L^2} \nn\\
&\leq 
C\|(f_h^{n+1})_{n=0}^m\|_{\ell^2(L^{q/2})} 
+C\|\psi_h^0\|_{H^1}  .
\end{align*} 
Again, the Sobolev embedding theorem implies
\begin{align}\label{thetaW1qa}
\begin{aligned}
\|(\theta^{n+1})_{n=0}^m\|_{\ell^2(W^{1,3/(1-\alpha)})}  
&\leq C\|(\theta^{n+1})_{n=0}^m\|_{\ell^2(H^{3/2+\alpha})}  \\
&\leq 
C\|(f_h^{n+1})_{n=0}^m\|_{\ell^2(L^{q/2})} 
+C\|\psi_h^0\|_{H^1} .
\end{aligned}
\end{align}

Comparing \eqref{EQDPSI} and \eqref{EQDPSI2},
we have 
\begin{align*}%\label{ErrPSI1}
\big(\eta D_\tau(\theta^{n+1}-\psi_h^{n+1}),
\varphi_h\big) 
+\frac{1}{\kappa^2}\big(\nabla(\theta^{n+1}-\psi_h^{n+1}),
\nabla\varphi_h\big) 
=0 ,\quad\forall\,\varphi_h\in {\mathbb S}_h^r ,
\end{align*}
which indicates that $\psi^{n+1}_h$
is the finite element approximation of $\theta^{n+1}$. 
The standard energy error estimate gives
\begin{align*} 
&\|(P_h\theta^{n+1}-\psi_h^{n+1})_{n=0}^m\|_{\ell^\infty(L^2)}
+\|(P_h\theta^{n+1}-\psi_h^{n+1})_{n=0}^m\|_{\ell^2(H^1)}  \\
&\leq C\|(P_h\theta^{n+1}-\theta^{n+1})_{n=0}^m\|_{\ell^2(H^1)}  \\
&
\leq C\|(\theta^{n+1})_{n=0}^m\|_{\ell^2(H^{3/2+\alpha})}h^{1/2+\alpha} \\
&\leq C(\|(f_h^{n+1})_{n=0}^m\|_{\ell^2(L^{q/2})} 
+\|\psi_h^0\|_{H^1})h^{1/2+\alpha} ,
\end{align*}
and by using the inverse inequality we derive 
\begin{align}\label{thetahW1qa}
\begin{aligned}
\|(P_h\theta^{n+1}-\psi_h^{n+1})_{n=0}^m\|_{\ell^2(W^{1,3/(1-\alpha)})} 
&\le Ch^{-1/2-\alpha}\|(P_h\theta^{n+1}-\psi_h^{n+1})_{n=0}^m\|_{\ell^2(H^1)} \\ 
&\le C(\|(f_h^{n+1})_{n=0}^m\|_{\ell^2(L^{q/2})} 
+\|\psi_h^0\|_{H^1}) .
\end{aligned}
\end{align} 
From \eqref{Defalpha} we know that 
$3/(1-\alpha)=q/(2-q/3)=q+\delta_q$ for some $\delta_q>0$. 
Since the $L^2$ projection operator 
$P_h$ is bounded on $W^{1,q+\delta_q}$,
the inequalities \eqref{thetaW1qa} and \eqref{thetahW1qa} 
imply \eqref{LemW1qD}.

The proof of Lemma \ref{LemmW1qD} is complete.  
\endproof\bigskip

We rewrite \eqref{FEM1} as
\begin{align}\label{ReWrEq}
\eta D_\tau\psi_h^{n+1}
-\frac{1}{\kappa^2}\Delta_h\psi^{n+1}_h
+\frac{i}{\kappa}P_h\big(\nabla\psi^{n+1}_h\cdot{\bf A}_h^{n+1}\big)
+\frac{i}{\kappa}\nabla_h\cdot\big(\psi^{n+1}_h {\bf A}_h^{n+1}\big)  &
  \nn\\
+P_h\Big( |\mathbf{A}^{n+1}_h|^2\psi^{n+1}_h + 
(|\psi^{n+1}_h|^{2}-1) \psi^{n+1}_h 
+i\eta\kappa \Theta(\psi_h^{n})\phi_h^n\Big) &= 0 ,
\end{align}
where the discretes operators 
\begin{align*}
&\Delta_h:{\mathbb S}_h^r\rightarrow {\mathbb S}_h^r, \\
&\nabla_h\cdot:{\cal L}^2\times {\cal L}^2
\times {\cal L}^2\rightarrow {\mathbb S}_h^r,\\
&P_h:{\cal L}^2\rightarrow {\mathbb S}_h^r
\end{align*}
are defined via duality by 
\begin{align*}
&(\Delta_hu_h,v_h)=-(\nabla u_h,\nabla v_h) , &&
\forall\, u_h,v_h\in {\mathbb S}_h^r ,\\
&(\nabla_h\cdot {\bf u},v_h)=-({\bf u}_h,\nabla v_h) , &&
\forall\, {\bf u}\in {\cal L}^2\times {\cal L}^2\times {\cal L}^2,
\,\, v_h\in {\mathbb S}_h^r ,\\
&(P_hu,v_h)=(u,v_h) , &&
\forall\, u\in {\cal L}^2,\,\, v_h\in {\mathbb S}_h^r .
\end{align*}
By applying Lemma \ref{LemmW1qD} to \eqref{ReWrEq}, 
using H\"older's inequality and \eqref{AhLinftLq}-\eqref{psihLinfH1}, 
we obtain
\begin{align}\label{psihn1a}
&\big\|\big(\psi^{n+1}_h\big)_{n=0}^m\big\|_{\ell^2(W^{1,q+\delta_q})} \nn\\
&\leq C\|\psi^{0}_h\|_{H^1}+
C\big\|\big(\nabla\psi^{n+1}_h\cdot{\bf A}_h^{n+1}\big)_{n=0}^m \big\|_{\ell^2(L^{q/2})} \nn\\
&\quad 
+C\big\|\big(\nabla_h\cdot(\psi^{n+1}_h {\bf A}_h^{n+1})\big)_{n=0}^m\big\|_{\ell^2(L^{q/2})}  \nn\\
&\quad
+C\big\|\big(|\mathbf{A}^{n+1}_h|^2\psi^{n+1}_h + 
(|\psi^{n+1}_h|^{2}-1) \psi^{n+1}_h 
-i\eta\kappa \Theta(\psi_h^{n})\phi^n \big)_{n=0}^m\big\|_{\ell^2(L^{q/2})} \nn\\
&\leq C+C\big\|\big(\nabla\psi^{n+1}_h\big)_{n=0}^m\big\|_{\ell^2(L^{q})}
\big\|\big({\bf A}_h^{n+1} \big)_{n=0}^m\big\|_{\ell^\infty(L^q)} \nn\\
&\quad 
+C\big\|\big(\nabla_h\cdot(\psi^{n+1}_h {\bf A}_h^{n+1})\big)_{n=0}^m\big\|_{\ell^2(L^{q/2})}  \nn\\
&\quad 
+C\big\|\big({\bf A}_h^{n+1}\big)_{n=0}^m \big\|_{\ell^\infty(L^q)}^2
\big\|\big(\psi^{n+1}_h\big)_{n=0}^m\big\|_{\ell^2(L^\infty)} \nn\\
&\quad  
+C\big(\big\|\big(\psi^{n+1}_h\big)_{n=0}^m\big\|_{\ell^{6}(L^{3q/2})}^3
\!+\!\big\|\big(\psi^{n+1}_h\big)_{n=0}^m\big\|_{\ell^2(L^{q/2})}   
\!+\!\big\|\big(\phi_h^n\big)_{n=0}^m\big\|_{\ell^2(L^{q/2})} \big) \nn\\
&\leq C+C\big(\big\|\big(\psi^{n+1}_h\big)_{n=0}^m\big\|_{\ell^2(W^{1,q})} 
+\big\|\big(\psi^{n+1}_h\big)_{n=0}^m\big\|_{\ell^2(L^\infty)} \big) \nn\\
&\quad 
+C\big\|\big(\nabla_h\cdot(\psi^{n+1}_h {\bf A}_h^{n+1})\big)_{n=0}^m\big\|_{\ell^2(L^{q/2})}  \nn\\
&\quad 
+C\big(\big\|\big(\psi^{n+1}_h\big)_{n=0}^m\big\|_{\ell^6(L^{3q/2})}^3
+\big\|\big(\psi^{n+1}_h\big)_{n=0}^m\big\|_{\ell^2(L^{q/2})} \big)\nn\\
&\leq C+\epsilon\big\|\big(\psi^{n+1}_h\big)_{n=0}^m\big\|_{\ell^2(W^{1,q+\delta_q})}
+C_\epsilon\big\|\big(\psi^{n+1}_h\big)_{n=0}^m\big\|_{\ell^2(H^1)} \nn\\
&\quad 
+C\big\|\big(\nabla_h\cdot(\psi^{n+1}_h {\bf A}_h^{n+1})\big)_{n=0}^m\big\|_{\ell^2(L^{q/2})}  \nn\\
&\quad 
 +C\big(\big\|\big(\psi^{n+1}_h\big)_{n=0}^m\big\|_{\ell^\infty(H^1)}^3
+\big\|\big(\psi^{n+1}_h\big)_{n=0}^m\big\|_{\ell^\infty(H^1)} \big)\nn\\
&\leq C_\epsilon
+\epsilon\big\|\big(\psi^{n+1}_h\big)_{n=0}^m\big\|_{\ell^2(W^{1,q+\delta_q})}
+C\big\|\big(\nabla_h\cdot(\psi^{n+1}_h {\bf A}_h^{n+1})\big)_{n=0}^m\big\|_{\ell^2(L^{q/2})}  \, ,   
\end{align}
%where $3q/2<6$ and $q^*<6$ satisfies 
%\begin{align}
%1/q^*+1/2=2/q ,
%\end{align}
where we have used the following interpolation inequality: 
\begin{align*}
\|(\psi^{n+1}_h)_{n=0}^m\|_{\ell^2(L^\infty)}
+\|(\psi^{n+1}_h)_{n=0}^m\|_{\ell^2(W^{1,q})}
%&\leq C\|\psi^{n+1}_h\|_{l^2_{n,m}(W^{1,q})}\\
&\leq \epsilon\|(\psi^{n+1}_h)_{n=0}^m\|_{\ell^2(W^{1,q+\delta_q})}
+C_\epsilon\|(\psi^{n+1}_h)_{n=0}^m\|_{\ell^2(H^1)} .
\end{align*}
To estimate $\|\nabla_h\cdot(\psi^{n+1}_h {\bf A}_h^{n+1})\|_{L^{q/2}}$ 
on the right-hand side of \eqref{psihn1a},
we let $q^*<6$ be the number satisfying 
$
1/q^*+1/2=2/q 
$
and use a duality argument:
for any $\eta_h\in{\mathbb S}_h^r$ we have 
\begin{align}\label{CondrkN}
&(\nabla_h\cdot(\psi^{n+1}_h {\bf A}_h^{n+1}),\eta_h) \nn\\
&=-(\psi^{n+1}_h {\bf A}_h^{n+1},\nabla\eta_h) \nn\\
&=({\bf A}_h^{n+1},\eta_h\nabla\psi^{n+1}_h )
-({\bf A}_h^{n+1},\nabla(\psi^{n+1}_h\eta_h) ) \nn\\
&=({\bf A}_h^{n+1},\eta_h\nabla\psi^{n+1}_h )
-(\phi_h^{n+1}, \psi^{n+1}_h\eta_h )
\qquad\qquad\qquad\,\, \mbox{by using \eqref{FEM3}
and \eqref{Condrk}} \nn\\
&\leq\|{\bf A}_h^{n+1}\|_{L^q} \|\nabla\psi^{n+1}_h\|_{L^q}
\|\eta_h\|_{L^{(q/2)'}} 
+\|\phi_h^{n+1}\|_{L^2}\|\psi^{n+1}_h\|_{L^{q^*}}
\|\eta_h\|_{L^{(q/2)'}} \nn\\
&\leq C\|\nabla\psi^{n+1}_h\|_{L^q}
\|\eta_h\|_{L^{(q/2)'}} 
+C\|\psi^{n+1}_h\|_{L^{q^*}}
\|\eta_h\|_{L^{(q/2)'}} ,
\qquad\qquad\mbox{by using \eqref{psihLinfH1}}
\end{align}
which implies 
\begin{align*} 
\|\nabla_h\cdot(\psi^{n+1}_h {\bf A}_h^{n+1}) \|_{L^{q/2}} 
&\leq C(\|\nabla\psi^{n+1}_h\|_{L^q} + \|\psi^{n+1}_h\|_{L^{q^*}})\\
&\leq C(\|\psi^{n+1}_h\|_{W^{1,q}} + \|\psi^{n+1}_h\|_{H^{1}}) ,
\end{align*}
and so 
\begin{align*}%\label{H1psiA}
\big\|\big(\nabla_h\cdot(\psi^{n+1}_h {\bf A}_h^{n+1})\big)_{n=0}^m\big\|_{\ell^2(L^{q/2})}
 &\leq 
C\big\|\big(\psi^{n+1}_h\big)_{n=0}^m\big\|_{\ell^2(W^{1,q})}
+C\big\|\big(\psi^{n+1}_h\big)_{n=0}^m\big\|_{\ell^2(H^1)}
 \nn\\
&\leq \epsilon\big\|\big(\psi^{n+1}_h\big)_{n=0}^m\big\|_{\ell^2(W^{1,q+\delta_q})}
+C_\epsilon\big\|\big(\psi^{n+1}_h\big)_{n=0}^m\big\|_{\ell^2(H^1)} \nn\\
&\leq \epsilon\big\|\big(\psi^{n+1}_h\big)_{n=0}^m\big\|_{\ell^2(W^{1,q+\delta_q})}
+C_\epsilon \qquad 
\mbox{by using \eqref{psihLinfH1}} ,
\end{align*} 
which together with \eqref{psihn1a} 
implies  
\begin{align}\label{DeltaPsi}
\big\|\big(\psi^{n+1}_h\big)_{n=0}^m\big\|_{\ell^2(W^{1,q+\delta_q})} 
\leq C  .
\end{align} 

For any $1\leq p\leq\infty$,
the space $\ell^p_{m}(W^{1,q})$
can be viewed as a subspace of $L^p(0,t_{m+1};W^{1,q})$
consisting of piecewise constant functions on
each subinterval $(t_n,t_{n+1}]$. Since 
$$
L^2(0,t_{m+1};W^{1,q})\cap L^\infty(0,t_{m+1};H^1)
\hookrightarrow L^{2/(1-\theta)}(0,t_{m+1};W^{1,q_\theta })
\quad\mbox{for any $\theta\in(0,1)$}, 
$$
with $
\frac{1}{q_\theta}=\frac{1-\theta}{q}+\frac{\theta}{2} $ 
(see \cite[page 106]{BL} on the 
complex interpolation of vector-valued $L^p$ spaces), 
it follows that
$\ell_{m}^2(W^{1,q})\cap \ell_{m}^\infty(H^1)
\hookrightarrow \ell_{m}^{2/(1-\theta)}(W^{1,q_\theta })$. 
By choosing $\theta$ to be sufficiently small we have
$3<q_\theta<q$ and so  
\begin{align*}
\big\|\big(\psi_h^{n+1}\big)_{n=0}^m\big\|_{\ell^{2/(1-\theta)}(L^{\infty})}
&\leq
C\big\|\big(\psi_h^{n+1}\big)_{n=0}^m\big\|_{\ell^{2/(1-\theta)}(W^{1,q_\theta})} \\
&\leq C\big\|\big(\psi_h^{n+1}\big)_{n=0}^m\big\|_{\ell^2(W^{1,q})} 
+C\big\|\big(\psi_h^{n+1}\big)_{n=0}^m\big\|_{\ell^\infty(H^1)} \\
&\leq C .
\end{align*}
In other words, we have 
\begin{align}
\sum_{n=0}^m\tau\|\psi_h^{n+1}\|_{L^{\infty}}^{2/(1-\theta)} 
\leq C_0
\quad\implies
\quad
\|\psi_h^{n+1}\|_{L^{\infty}}
\le (C_0^{(1-\theta)/2}\tau^{\theta/2}) \tau^{-1/2}
\end{align} 
for some positive constant $C_0$
(which is independent of $m$). 
When $\tau<\tau_0:=C_0^{-(1-\theta)/\theta}$, 
we have $C_0^{(1-\theta)/2}\tau^{\theta/2}<1$ and 
the last inequality implies 
\eqref{MathInd} for $n=m+1$. 
Hence, the mathematical induction on 
\eqref{MathInd} is completed under the condition
$\tau<\tau_0$.  
As a consequence, \eqref{MathInd}-\eqref{psihLinfH1} and \eqref{DeltaPsi} hold 
for $m=N-1$.

Substituting ${\bf a}_h=\nabla\phi_h^{n+1}$
in \eqref{FEM2} 
and using \eqref{FEM3}, we obtain 
\begin{align}
\begin{aligned}
(D_{\tau}\phi_{h}^{n+1},\phi_{h}^{n+1}) 
+\frac{1}{2}\|\nabla\phi_{h}^{n+1}\|_{L^2}^2
&\leq
C\biggl\|{\rm Re}\bigg[\overline\psi_h^{n}\bigg(\frac{i}{\kappa} \nabla 
+{\bf A}_h^{n}\bigg) \psi_h^{n}\bigg]\bigg\|_{L^2}^2 \\
%&\leq C\|\psi_h^{n+1}\|_{L^6}^2 
%\bigg\|\frac{i}{\kappa} \nabla \psi_h^{n+1}
%+{\bf A}_h^{n}\psi_h^{n+1}\bigg\|_{L^3}^2 \nn\\
&\leq C\|\psi_h^{n}\|_{L^6}^2 
(\|\nabla \psi_h^{n}\|_{L^3}^2
+\|{\bf A}_h^{n}\|_{L^3}^2\|\psi_h^{n}\|_{L^\infty}^2) \\
&\leq C\|\psi_h^{n}\|_{H^1}^2 
(\|\psi_h^{n}\|_{W^{1,3}}^2
+C\|\psi_h^{n}\|_{W^{1,q}}^2) \\
&\leq C\|\psi_h^{n}\|_{H^1}^2 \|\psi_h^{n}\|_{W^{1,q}}^2 .
\end{aligned}
\end{align}
Summing up the inequality above for $n=0, 1,\dots,N-1$, 
and using \eqref{DeltaPsi} with $m=N-1$, we obtain 
\begin{align}\label{H1phi}
\|(\nabla\phi_{h}^{n+1})_{n=0}^{N-1}\|_{\ell^2(L^2)}^2
&\leq C\|(\psi_h^{n})_{n=0}^{N-1}\|_{\ell^\infty(H^1)}^2 
\|(\psi_h^{n})_{n=0}^{N-1}\|_{\ell^2(W^{1,q})}^2
\leq C .
\end{align}
Then substituting ${\bf a}_h=\nabla\chi_h$
in \eqref{FEM2}, we obtain
\begin{align}
&(D_{\tau}\phi_{h}^{n+1},\chi_{h}) 
+(\nabla\phi_{h}^{n+1} \, , \nabla\chi_{h} )
+ {\rm Re}\bigg(\overline\psi_h^{n}\bigg(\frac{i}{\kappa} \nabla 
+{\bf A}_h^{n}\bigg) \psi_h^{n} , \nabla\chi_{h}\bigg) =0  \, , 
\end{align}
which implies  
\begin{align}\label{H1phi2}
\begin{aligned}
&\big\|\big(D_\tau\phi_{h}^{n+1}\big)_{n=0}^{N-1}\big\|_{\ell^2((H^1)')}  \\
&\leq C
\bigg(
\big\|\big(\nabla\phi_h^{n+1}\big)_{n=0}^{N-1}\big\|_{\ell^2(L^2)}
+\biggl\|\bigg({\rm Re}\bigg[\overline\psi_h^{n}\bigg(\frac{i}{\kappa} \nabla 
+{\bf A}_h^{n}\bigg) \psi_h^{n}\bigg]\bigg)_{n=0}^{N-1}\bigg\|_{\ell^2(L^2)}  \bigg)
\leq C . 
\end{aligned} 
\end{align} 
via duality.  
The proof of Lemma \ref{UniFEst} is complete. \endproof

\subsection{Compactness
of the finite element solution}\label{CompFES}

For $t\in[t_n,t_{n+1}]$, $n=0,1,\cdots,N-1$,  
we define 
\begin{align*}
&\psi_{h,\tau}(t)
=\frac{1}{\tau}[(t_{n+1}-t)\psi_h^n+(t-t_n)\psi_h^{n+1}] ,\\
&{\bf A}_{h,\tau}(t)=\frac{1}{\tau}[(t_{n+1}-t){\bf A}_h^n
+(t-t_n){\bf A}_h^{n+1}] ,\\\
&\phi_{h,\tau}(t)=\frac{1}{\tau}[(t_{n+1}-t)
\phi_h^n+(t-t_n)\phi_h^{n+1}] .
%&{\bf F}_{h,\tau}(t)= 
%{\rm Re}\bigg[\overline\psi_{h,\tau}(t)
%\bigg(\frac{i}{\kappa}\nabla\psi_{h,\tau}(t)
%+ {\bf A}_{h,\tau}(t)\psi_{h,\tau}(t)\bigg)\bigg] ,\\
%%&\psi_{h,\tau}^+(t)=\psi_h^{n+1},\\
%%&\psi_{h,\tau}^-(t)=\psi_h^{n},\\
%%&{\bf A}_{h,\tau}^+(t)={\bf A}_h^{n+1},\\
%%&{\bf A}_{h,\tau}^-(t)={\bf A}_h^{n} ,\\
%&{\bf F}_{h,\tau}^*(t)= {\rm Re}
%\bigg[\psi_{h}^{n+1}\bigg(\frac{i}{\kappa}\nabla\psi_{h}^{n+1}
%+ {\bf A}_{h}^n\psi_{h}^{n+1}\bigg)\bigg] .
\end{align*} 
In other words,
$\psi_{h,\tau}$, ${\bf A}_{h,\tau}$
and ${\bf B}_{h,\tau}$ are the piecewise linear interpolation
of the functions $\psi_h^n$,
${\bf A}_h^n$ and ${\bf B}_h^n$ 
on the interval $[0,T]$, 
respectively. 
Then \eqref{LUniFEst} implies  
\begin{align}\label{Epsi3} 
&\|\psi_{h,\tau}\|_{H^1(0,T;L^2)}  
+\|\psi_{h,\tau}\|_{L^\infty(0,T;H^1)}  
+ \|\psi_{h,\tau}\|_{L^2(0,T;L^\infty)}
+ \|\psi_{h,\tau}\|_{L^2(0,T;W^{1,q})}
\leq C ,\\
&\|{\bf A}_{h,\tau}\|_{H^1(0,T;L^2)}
+\|{\bf A}_{h,\tau}\|_{L^\infty(0,T;L^q)}   
+\|\nabla\times{\bf A}_{h,\tau}\|_{L^\infty(0,T;L^2)}    \leq  C ,
\label{Epsi3-2} \\
&\|\phi_{h,\tau}\|_{L^\infty(0,T;L^2)}  
+\|\phi_{h,\tau}\|_{L^2(0,T;H^1)}
+\|\partial_t\phi_{h,\tau}\|_{L^2(0,T;(H^1)')}  \leq C .
\label{Epsi3-3}
\end{align} 

We see that $\psi_{h,\tau}$
is bounded in $L^\infty(0,T;{\cal H}^1)\cap H^{1}(0,T;{\cal L}^2)
\hookrightarrow C^{\theta/2}([0,T];{\cal H}^{1-\theta})$
for any $\theta\in(0,1)$. 
Since for any given $1<p<6$ there is a small $\theta$ 
such that $C^{\theta/2}([0,T];{\cal H}^{1-\theta})$ 
is compactly embedded 
into $C([0,T];{\cal L}^p)$, 
\eqref{Epsi3} implies compactness of $\psi_{h,\tau}$
in $C([0,T];{\cal L}^p)$ for any $1<p<6$. 
Hence, for any sequence $(h_m,\tau_m)\rightarrow (0,0)$, 
the inequality \eqref{Epsi3} implies the existence of 
a subsequence, also denoted by $(h_m,\tau_m)$
for the simplicity of the notations, which satisfies  
\begin{align}
&\partial_t\psi_{h_m,\tau_m}\rightarrow \partial_t\Psi
&&\mbox{weakly in $L^2(0,T;{\cal L}^2)$}, \label{Convpsi1}\\
&\psi_{h_m,\tau_m}\rightarrow \Psi 
&&\mbox{weakly$^*$ in $L^\infty(0,T;{\cal H}^1)$ }, \label{Convpsi1-2}\\
&\psi_{h_m,\tau_m}\rightarrow \Psi 
&&\mbox{weakly in $L^2(0,T;{\cal W}^{1,q})$
for some $q>3$}, \label{Convpsi2} \\ 
&\psi_{h_m,\tau_m}\rightarrow \Psi 
&&\mbox{strongly in $C([0,T];{\cal L}^p)$
for any $1<p<6$} . 
\label{Convpsi3}
\end{align}
for some function $\Psi$. 

Using the notation of Definition \ref{DisDiv}, 
we have $\phi_{h,\tau}=\nabla_h^{\mathbb N}\cdot{\bf A}_{h,\tau}$ 
and \eqref{Epsi3-2}-\eqref{Epsi3-3} 
imply that ${\bf A}_{h,\tau}$ is bounded in 
the norm of 
$$L^\infty(0,T;{\bf H}_h({\rm curl,div}))
\cap H^1(0,T;{\bf L}^2)
\hookrightarrow 
C^{\theta/2}([0,T];{\bf Y}_{1-\theta}) ,
\quad\forall\,\theta\in(0,1) ,
$$
where 
${\bf Y}_{1-\theta}:=({\bf H}_h({\rm curl,div}),{\bf L}^2)_{1-\theta}$
is the real interpolation space between
${\bf H}_h({\rm curl,div})$ and ${\bf L}^2$
(see \cite{BL}). 
Lemma \ref{SobolevD} says that 
a set of functions which are bounded in the norm 
of ${\bf H}_h({\rm curl,div})$
is compact in ${\bf L}^2$, which implies that 
a set of functions which are bounded in the norm 
of the interpolation space 
${\bf Y}_{1-\theta}$
is also compact in ${\bf L}^2$
(see Theorem 3.8.1, page 56 of \cite{BL}).  
Hence, $C^{\theta/2}([0,T];{\bf Y}_{1-\theta})$
is compactly embedded into $C([0,T];{\bf L}^2)$, 
and for any sequence ${\bf A}_{h_m,\tau_m}$
there exists a subsequence which converges
to some function 
${\bf\Lambda}$ strongly in $C([0,T];{\bf L}^2)$.
On the other hand, since 
${\bf H}_h({\rm curl,div})\hookrightarrow {\bf L}^{q+\delta}$
for some $q>3$ and $\delta>0$, 
by choosing $\theta$ small enough we have
$C^{\theta/2}([0,T];{\bf Y}_{1-\theta})
\hookrightarrow 
C([0,T];{\bf L}^{q+\delta/2})$. 
The boundedness of ${\bf A}_{h,\tau}$ in 
$C([0,T];{\bf L}^{q+\delta/2})$ implies the existence of 
a subsequence of ${\bf A}_{h_m,\tau_m}$ which converges
weakly$^*$ to some function in 
$L^\infty(0,T;{\bf L}^{q+\delta/2})$. 
This weak limit must also be ${\bf \Lambda}$,
and 
\begin{align}
\|{\bf A}_{h_m,\tau_m}-{\bf \Lambda}\|_{L^\infty(0,T;{\bf L}^{q})}
&\leq \|{\bf A}_{h_m,\tau_m}-{\bf \Lambda}\|_{L^\infty(0,T;{\bf L}^2)}^{1-\theta}
\|{\bf A}_{h_m,\tau_m}-{\bf \Lambda}\|_{L^\infty(0,T;{\bf L}^{q+\delta/2})}^\theta \nn\\
&\leq C\|{\bf A}_{h_m,\tau_m}-{\bf \Lambda}\|_{L^\infty(0,T;{\bf L}^2)}^{1-\theta}
\end{align}
for some $\theta>0$.
In other words, 
${\bf A}_{h_m,\tau_m}\in C([0,T];{\bf L}^{q})$ converges
to  ${\bf\Lambda}$ strongly in $L^\infty(0,T;{\bf L}^q)$,
which implies  
${\bf\Lambda}\in C([0,T];{\bf L}^{q})$.
To conclude, there exists a subsequence of $(h_m,\tau_m)$,
which is also denoted by $(h_m,\tau_m)$
for the simplicity of the notations, such that 
\begin{align}
&\partial_t{\bf A}_{h_m,\tau_m}\rightarrow \partial_t{\bf \Lambda}
&&\mbox{weakly in $L^2(0,T;{\bf L}^2)$}, \label{ConvA1}\\
&\nabla\times{\bf A}_{h_m,\tau_m}\rightarrow \nabla\times {\bf\Lambda}
&&\mbox{weakly$^*$ in $L^\infty(0,T;{\bf L}^2)$}, \label{ConvA3}\\ 
&{\bf A}_{h_m,\tau_m}\rightarrow {\bf \Lambda}
&&\mbox{strongly in $C([0,T];{\bf L}^q)$
for some $q>3$} , 
\label{ConvA4}
\end{align}
for some function ${\bf\Lambda}$. 

Similarly, \eqref{Epsi3-3} implies 
the existence of a subsequence
such that  
\begin{align}
&\phi_{h_m,\tau_m}\rightarrow \Phi
&&\mbox{weakly$^*$ in $L^\infty(0,T;L^2)$}, \label{ConvA2}\\ 
&\phi_{h_m,\tau_m}\rightarrow \Phi
&&\mbox{weakly in $L^2(0,T;H^1)$}, \label{ConvA2-2}\\ 
&\phi_{h_m,\tau_m}\rightarrow \Phi
&&\mbox{strongly in $L^2(0,T;L^2)$} . \label{ConvA2-2-2} 
\end{align}
for some function $\Phi$. 

For any $\chi\in L^2(0,T;H^1)$ 
and finite element functions 
$\chi_{h_m,\tau_m}\rightarrow \chi$ in $L^2(0,T;H^1)$, 
equation \eqref{FEM3} implies 
\begin{align}
&\int_0^T (\phi_{h_m,\tau_m},\chi) \d t=
\int_0^T\bigg[(\phi_{h_m,\tau_m},\chi-\chi_{h_m,\tau_m})
+({\bf A}_{h_m,\tau_m},\nabla\chi_{h_m,\tau_m})\bigg]\d t
\end{align}
As $h_m,\tau_m\rightarrow 0$,
the equation above tends to
\begin{align}
&\int_0^T(\Phi,\chi)\d t=
\int_0^T({\bf \Lambda},\nabla\chi) \d t ,
\end{align}
which implies that 
\begin{align}
\nabla\cdot{\bf \Lambda}=-\Phi\in L^\infty(0,T;L^2)\cap 
L^2(0,T;H^1) .
\end{align}

Now we consider compactness of $\psi_{h,\tau}^\pm$, 
${\bf A}_{h,\tau}^\pm$ and $\phi_{h,\tau}^\pm$
by utilizing the compactness of 
$\psi_{h,\tau}$, 
${\bf A}_{h,\tau}$ and $\phi_{h,\tau}$. 
Since $\psi_{h,\tau}$ is bounded in $H^1(0,T;L^2)\cap L^\infty(0,T;H^1)\hookrightarrow C^{(1-\theta)/2}([0,T];L^{p_\theta})$ for 
$$
\frac{1}{p_\theta}=\frac{1-\theta}{2}+\frac{\theta}{6},\qquad
\forall\,\theta\in(0,1), 
$$ 
it follows that 
\begin{align}
\|\psi_{h,\tau}(t)-\psi_{h,\tau}^+(t)\|_{L^{p_\theta}}
&=\bigg\|\frac{t_{n+1}-t}{\tau} 
(\psi_{h,\tau}(t_n)-\psi_{h,\tau}(t_{n+1}))\bigg\|_{L^{p_\theta}}\nn\\[5pt]
&\leq C\|\psi_{h,\tau}\|_{C^{(1-\theta)/2}([0,T];L^{p_\theta})} \tau^{(1-\theta)/2}  
\end{align} 
for $t\in(t_n,t_{n+1})$, and so
\begin{align}
\|\psi_{h,\tau}-\psi_{h,\tau}^+\|_{L^\infty(0,T;L^{p_\theta})}
\leq C\|\psi_{h,\tau}\|_{C^{(1-\theta)/2}([0,T];L^{p_\theta})} \tau^{(1-\theta)/2}
\rightarrow 0\quad\mbox{as}\,\,\, \tau\rightarrow 0 . 
\end{align} 
Similarly, we also have
\begin{align}
\|\psi_{h,\tau}-\psi_{h,\tau}^-\|_{L^\infty(0,T;L^{p_\theta})}
\leq C\|\psi_{h,\tau}\|_{C^{\alpha_p}([0,T];L^{p_\theta})} \tau^{(1-\theta)/2}
\rightarrow 0\quad\mbox{as}\,\,\, \tau\rightarrow 0 . 
\end{align} 
Since $\psi_{h_m,\tau_m}$ converges
strongly in $L^\infty(0,T;L^{p_\theta})$,
it follows that both $\psi_{h_m,\tau_m}^-$
and $\psi_{h_m,\tau_m}^+$ converge to the same function 
strongly in $L^\infty(0,T;L^{p_\theta})$. 
Hence, there exists a subsequence 
which satisfies 
\begin{align}
&\psi_{h_m,\tau_m}^{\pm}\rightarrow \Psi 
&&\mbox{weakly$^*$ in $L^\infty(0,T;H^1)$ }, \label{Convpsi1-3}\\ 
&\psi_{h_m,\tau_m}^{\pm}\rightarrow \Psi 
&&\mbox{weakly in $L^2(0,T;W^{1,q})$
for some $q>3$} , \label{Convpsi3-3}\\
&\psi_{h_m,\tau_m}^{\pm}\rightarrow \Psi 
&&\mbox{strongly in $L^\infty(0,T;L^p)$
for any $1<p<6$} \label{Convpsi2-3} .
\end{align}
In a similar way one can prove
\begin{align}
&{\bf A}_{h_m,\tau_m}^{\pm}\rightarrow {\bf \Lambda}
&&\mbox{strongly in $L^\infty(0,T;L^q)$
for some $q>3$} , \label{ConvA4-3} \\
&\nabla\times{\bf A}_{h_m,\tau_m}^{\pm}\rightarrow
\nabla\times{\bf \Lambda}
&&\mbox{weakly$^*$ in $L^\infty(0,T;L^2)$}, \label{ConvA3-3}\\  
&\phi_{h_m,\tau_m}^{\pm}\rightarrow \Phi=-\nabla\cdot{\bf A}
&&\mbox{weakly$^*$ in $L^\infty(0,T;L^2)$}, \label{ConvA2-3}\\ 
&\phi_{h_m,\tau_m}^{\pm}\rightarrow \Phi
&&\mbox{weakly in $L^2(0,T;H^1)$} . \label{ConvA2-3-3}\\
&\phi_{h_m,\tau_m}^{\pm}\rightarrow \Phi
&&\mbox{strongly in $L^2(0,T;L^2)$} . \label{ConvA2-3-4}
\end{align} 
From \eqref{Convpsi1-3}-\eqref{ConvA4-3} and \eqref{ConvA2-3-3}
we see that 
\begin{align}\label{ConvFhm}  
&\psi_{h_m,\tau_m}^+\bigg(\frac{i}{\kappa}\nabla
+ {\bf A}_{h_m,\tau_m}^+\bigg)\psi_{h_m,\tau_m}^+\rightarrow 
\overline\Psi\bigg(\frac{i}{\kappa}\nabla + {\bf\Lambda}\bigg)\Psi  
&&\mbox{weakly in $L^2(0,T;L^2)$} ,\\
&\bigg(\frac{i}{\kappa}\nabla+{\bf A}_{h_m,\tau_m}^+\bigg)\psi_{h_m,\tau_m}^+ \rightarrow 
\bigg(\frac{i}{\kappa}\nabla + {\bf\Lambda}\bigg)\Psi
&&\mbox{weakly in $L^2(0,T;L^3)$},\\
&{\bf A}_{h_m,\tau_m}^+\cdot
\bigg(\frac{i}{\kappa}\nabla+{\bf A}_{h_m,\tau_m}^+\bigg)\psi_{h_m,\tau_m}^+ \rightarrow 
{\bf \Lambda}\cdot\bigg(\frac{i}{\kappa}\nabla+ {\bf\Lambda} \bigg)\Psi 
&&\mbox{weakly in $L^2(0,T;L^{3/2})$} ,\\
&\Theta(\psi_{h_m,\tau_m}^-)\phi_{h_m,\tau_m}^- 
\rightarrow \Theta(\Psi)\Phi
&&\mbox{weakly in $L^2(0,T;L^2)$},\\
&|\psi_{h_m,\tau_m}^+|^3 
\rightarrow |\Psi|^3 
&&\mbox{weakly in $L^2(0,T;L^2)$} .
\label{WWWA}
\end{align} 
Moreover, from \eqref{Convpsi3} and \eqref{ConvA4} we know that
$\Psi(\cdot,0)=\psi_0$ and ${\bf \Lambda}(\cdot,0)={\bf A}_0$.

\subsection{Convergence
to the PDE's solution}\label{ConvgS}

It remains to prove 
\begin{align}
\Psi=\psi,\qquad
{\bf\Lambda}={\bf A}\qquad\mbox{and}\qquad 
\Phi=\phi ,
\end{align} 
so that \eqref{Convpsi1-3}-\eqref{ConvA2-3-4} imply 
Theorem \ref{MainTHM}. 

For any given $\varphi\in L^2(0,T;{\cal H}^1)$,
we choose finite element functions 
$\varphi_{h,\tau}\in L^2(0,T;{\mathbb S}_h^r)$ which 
converge to $\varphi$ strongly 
in $L^2(0,T;{\cal H}^1)$ 
as $h\rightarrow 0$. 
Then \eqref{FEM1} implies  
\begin{align}
&\int_0^T\bigg[(\eta \partial_t\psi_{h,\tau}, \varphi_{h,\tau}) 
+ (i\eta \kappa
\Theta(\psi_{h,\tau}^-)\phi_{h,\tau}^- ,\varphi_{h,\tau})\bigg]\d t \nn \\
&
+\int_0^T\bigg[ \bigg( \bigg(\frac{i}{\kappa}\nabla+{\bf A}_{h,\tau}^+\bigg)\psi_{h,\tau}^+ \,, 
\bigg(\frac{i}{\kappa}\nabla+{\bf A}_{h,\tau}^+\bigg)\varphi_{h,\tau}\bigg)  
+((|\psi_{h,\tau}^+|^{2}-1)\psi_{h,\tau}^+,\varphi_{h,\tau})\bigg]\d t =0   . \nn
\end{align}
Let $h=h_m\rightarrow 0$ and 
$\tau=\tau_m\rightarrow 0$ in the equation above
and use \eqref{Convpsi1} 
and \eqref{Convpsi1-3}-\eqref{WWWA}. 
We obtain 
\begin{align}\label{PDEPsiF}
&\int_0^T\bigg[(\eta \partial_t\Psi, \varphi)  
+ (i\eta \kappa\Theta(\Psi)\Phi,\varphi) 
+ \bigg( \bigg(\frac{i}{\kappa}\nabla+{\bf \Lambda}\bigg)\Psi \,, 
\bigg(\frac{i}{\kappa}\nabla+{\bf \Lambda}\bigg)\varphi\bigg)  
\bigg]\d t\nn\\
&
+\int_0^T((|\Psi|^{2}-1)\Psi,\varphi) \d t=0, 
\end{align}
for any given $\varphi\in L^2(0,T;{\cal H}^1)$. 
Now we prove $|\Psi|\leq 1$ by using
the following lemma. 
\begin{lemma}\label{UnBDPsi}
{\it For any given ${\bf \Lambda}\in 
L^\infty(0,T;{\bf H}({\rm curl,div}))$
and $\Phi\in L^\infty(0,T;L^2)$,  
the nonlinear equation {\rm\eqref{PDEPsiF}}
has a unique weak solution 
$\Psi\in L^2(0,T;{\cal H}^1)
\cap H^1(0,T;({\cal H}^1)')$
under the initial condition 
$\Psi(\cdot,0)=\psi_0$.
Moreover, the solution  
satisfies that $|\Psi|\leq 1$ a.e. in $\Omega\times(0,T)$.
}
\end{lemma}
\noindent{\it Proof of Lemma \ref{UnBDPsi}}.$\quad$
%Lemma 3.2
%of \cite{LY} implies the existence of a weak solution
%$\Psi\in L^2(0,T;{\cal H}^1)
%\cap H^1(0,T;({\cal H}^1)')$
%satisfying $|\Psi|\leq 1$ a.e. in $\Omega\times(0,T)$. 
%It remains to prove that
%there exists at most one 
%weak solution of \eqref{PDEPsiF}
%with the regularity $\Psi\in L^2(0,T;{\cal H}^1)
%\cap H^1(0,T;({\cal H}^1)')$. 
To prove uniqueness of the solution, 
let us suppose that there are two solutions
$\Psi,\widetilde\Psi\in L^2(0,T;{\cal H}^{1})
\cap H^1(0,T;({\cal H}^1)')$ 
for the equation \eqref{PDEPsiF}
with the same initial condition. 
Then ${\cal E}=\Psi-\widetilde\Psi$
satisfies the equation 
\begin{align*} 
&\int_0^T(\eta \partial_t{\cal E}, \varphi) \d t
+\int_0^T (i\eta \kappa(\Theta(\Psi)-\Theta(\widetilde\Psi))\Phi,\varphi) 
\d t\nn\\
&
+ \int_0^T\bigg( \bigg(\frac{i}{\kappa}\nabla+{\bf \Lambda}\bigg){\cal E} \,, 
\bigg(\frac{i}{\kappa}\nabla+{\bf \Lambda}\bigg)\varphi\bigg)  \d t
+\int_0^T(|\Psi|^{2}\Psi-|\widetilde\Psi|^{2}\widetilde\Psi,\varphi) 
\d t=
\int_0^T({\cal E},\varphi) \d t
\end{align*}
for any $\varphi\in L^2(0,T;{\cal H}^1)$. 
Since 
\begin{align*}
|\Theta(\Psi)-\Theta(\widetilde\Psi)|\leq |{\cal E}|
\qquad
\mbox{and}\qquad
(|\Psi|^{2}\Psi-|\widetilde\Psi|^{2}\widetilde\Psi,
\Psi-\widetilde\Psi)\geq 0 ,
\end{align*} 
by substituting $\varphi(x,t)={\cal E}(x,t)1_{[0,s]}(t)$ 
into the equation above,
we obtain
\begin{align}\label{CalE1}
&\frac{\eta}{2} \|{\cal E}(\cdot,s)\|_{L^2}^2 
+ \int_0^s\bigg\|\bigg(\frac{i}{\kappa}\nabla+{\bf \Lambda}\bigg){\cal E}\bigg\|_{L^2}^2\d t  \nn\\
&\leq 
\int_0^s\|{\cal E}(\cdot,t)\|_{L^2}^2\d t
+C\|\Phi\|_{L^\infty(0,s;L^2)}\||{\cal E}|^2\|_{L^1(0,s;L^2)} \nn\\
&\leq 
\int_0^s\|{\cal E}(\cdot,t)\|_{L^2}^2\d t
+C\|{\cal E}\|_{L^2(0,s;L^4)}^2 \nn\\
&\leq 
C_\epsilon\int_0^s\|{\cal E}(\cdot,t)\|_{L^2}^2\d t
+\epsilon\int_0^s\|\nabla {\cal E}(\cdot,t)\|_{L^2}^2 \d t ,
\end{align}
where $\epsilon\in(0,1)$ is arbitrary.

Note that ${\bf\Lambda}\in L^\infty(0,T;{\bf H}({\rm curl,div}))
\hookrightarrow L^\infty(0,T;{\bf L}^q)$
for some $q>3$. 
If we let $\bar q<6$ be the number satisfying
$1/q+1/\bar q=1/2$
and let $\theta_q\in(0,1)$
be the number satisfying 
$1/q=(1-\theta_q)/2+\theta_q/6$, then 
\begin{align*}
\|\nabla{\cal E}\|_{L^2}
&\leq \kappa\bigg\|\frac{i}{\kappa}\nabla{\cal E}
+ {\bf \Lambda}{\cal E}\bigg\|_{L^2}
+\kappa \|{\bf \Lambda}{\cal E}\|_{L^2}  \\
&\leq \kappa\bigg\|\frac{i}{\kappa}\nabla{\cal E}
+ {\bf \Lambda}{\cal E}\bigg\|_{L^2}
+\kappa \|{\bf \Lambda}\|_{L^q}\|{\cal E}\|_{L^{\bar q}}  \\
&\leq \kappa\bigg\|\frac{i}{\kappa}\nabla{\cal E}
+ {\bf \Lambda}{\cal E}\bigg\|_{L^2}
+C\|{\cal E}\|_{L^2}^{1-\theta_q}\|{\cal E}\|_{L^6}^{\theta_q} \\
&\leq 
\kappa\bigg\|\frac{i}{\kappa}\nabla{\cal E}
+ {\bf\Lambda}{\cal E}\bigg\|_{L^2}
+ \epsilon\|\nabla{\cal E}\|_{L^2} 
+C_\epsilon\|{\cal E}\|_{L^2} ,
\end{align*}
which implies
\begin{align*} 
\frac{1}{2\kappa}\|\nabla{\cal E}\|_{L^2}
&\leq 
\bigg\|\frac{i}{\kappa}\nabla{\cal E}
+ {\bf\Lambda}{\cal E}\bigg\|_{L^2}
+C\|{\cal E}\|_{L^2} .
\end{align*}
Substituting the last inequality into \eqref{CalE1},
we obtain
\begin{align*} 
&\frac{\eta}{2} \|{\cal E}(\cdot,s)\|_{L^2}^2 
+ \frac{1}{2\kappa}\int_0^s\|\nabla{\cal E}(\cdot,t)\|_{L^2}^2\d t  
\leq 
C_\epsilon\int_0^s\|{\cal E}(\cdot,t)\|_{L^2}^2\d t
+\epsilon\int_0^s\|\nabla {\cal E}(\cdot,t)\|_{L^2}^2 \d t ,
\end{align*}
which further reduces to
(by choosing sufficiently small $\epsilon$) 
\begin{align*} 
&\frac{\eta}{2} \|{\cal E}(\cdot,s)\|_{L^2}^2 
+ \frac{1}{2\kappa}\int_0^s\|\nabla{\cal E}(\cdot,t)\|_{L^2}^2\d t  
\leq 
C\int_0^s\|{\cal E}(\cdot,t)\|_{L^2}^2\d t .
\end{align*}
By applying Gronwall's inequality we derive 
\begin{align*} 
&\max_{0\leq t\leq T}\|{\cal E}(\cdot,t)\|_{L^2}^2 
\leq 
C\|{\cal E}(\cdot,0)\|_{L^2}^2=0 ,
\end{align*}
which implies the uniqueness of the 
weak solution of \eqref{PDEPsiF}. \medskip

Under the regularity of ${\bf \Lambda}$
and $\Phi$, existence of weak solutions of 
the weak formulated equation
\begin{align}\label{PDEPsiF333}
&\int_0^T\bigg[(\eta \partial_t\Psi, \varphi)  
+ (i\eta \kappa \Psi \Phi,\varphi) 
+ \bigg( \bigg(\frac{i}{\kappa}\nabla+{\bf \Lambda}\bigg)\Psi \,, 
\bigg(\frac{i}{\kappa}\nabla+{\bf \Lambda}\bigg)\varphi\bigg)  
\bigg]\d t\nn\\
&
+\int_0^T((|\Psi|^{2}-1)\Psi,\varphi) \d t=0, 
\qquad\forall\,\varphi\in L^2(0,T;{\cal H}^1), 
\end{align}
is obvious if
one can prove the a priori estimate 
\begin{align}\label{PWEst111}
\mbox{$|\Psi|\leq 1$\,\, a.e.\, in\,\, $\Omega\times(0,T)$.}
\end{align}
To prove the above inequality, 
we let $(|\Psi|^2-1)_+$ denote the positive part of $|\Psi|^2-1$
and integrate this equation against $\overline\Psi(|\Psi|^2-1)_+$. 
By considering the real part of the result, for any $t'\in(0,T)$ we have 
\begin{align*}
& \int_\Omega
\bigg(\frac{\eta}{4}\big(|\Psi(x,t')|^2-1\big)_+ ^2\bigg)\d x
 + \int_0^{t'}\int_\Omega (|\Psi|^{2}-1)^2_+ |\Psi| ^2\d x\d t\\
&=-\int_0^{t'}{\rm Re}\int_\Omega 
\bigg(\frac{i}{\kappa} \nabla  \Psi+ {\bf\Lambda} \Psi\bigg)
\bigg(-\frac{i}{\kappa} \nabla 
+ {\bf\Lambda}\bigg)[\overline\Psi (|\Psi|^2-1)_+]\d x\d t\\
&=-\int_0^{t'}\int_\Omega \bigg|\frac{i}{\kappa} 
\nabla  \Psi+ {\bf\Lambda} \Psi\bigg|^2
 (|\Psi|^2-1)_+ \d x \d t\\
&\quad + \int_0^{t'}{\rm Re}\int_{\{|\Psi|^2>1\}} 
\bigg(\frac{i}{\kappa} \nabla  \Psi
+ {\bf\Lambda} \Psi\bigg)\overline\Psi \bigg(\frac{i}{\kappa}  
\Psi\nabla\overline\Psi +\frac{i}{\kappa}\overline\Psi \nabla\Psi \bigg)\d x\d t\\
&=-\int_0^{t'}\int_\Omega \bigg|\frac{i}{\kappa} 
\nabla  \Psi+ {\bf\Lambda} \Psi\bigg|^2
 (|\Psi|^2-1)_+ \d x\d t\\
&\quad -\int_0^{t'}{\rm Re}\int_{\{|\Psi|^2>1\}}(|\Psi|^2|\nabla\Psi|^2
+ (\overline\Psi )^2\nabla\Psi\cdot \nabla\Psi)\d x\d t\\
& \leq 0,
\end{align*}
which implies that $\int_\Omega(|\Psi(x,t')|^2-1)_+ ^2 \d x
=0$, and this gives \eqref{PWEst111}. 
Since $|\Psi|\leq 1$, it follows that 
$\Theta(\Psi)=\Psi$ and so \eqref{PDEPsiF333}
reduces to \eqref{PDEPsiF}. This proves the existence
of weak solutions for \eqref{PDEPsiF} satisfying
$|\Psi|\leq 1$. 

The proof of Lemma \ref{UnBDPsi} is complete.
\endproof\bigskip

Lemma \ref{UnBDPsi} implies
\begin{align}
|\Psi|\leq 1 \quad\mbox{a.e. in $\Omega\times(0,T)$, }
\end{align}
which together with \eqref{PDEPsiF} implies 
\begin{align}\label{PDEPsiFM}
&\int_0^T\bigg[ (\eta \partial_t\Psi, \varphi) 
+ (i\eta \kappa \Psi\Phi,\varphi) 
+ \bigg( \bigg(\frac{i}{\kappa}\nabla+{\bf \Lambda}\bigg)\Psi \,, 
\bigg(\frac{i}{\kappa}\nabla+{\bf \Lambda}\bigg)\varphi\bigg)\bigg]\d t  \nn\\
&
+\int_0^T ((|\Psi|^{2}-1)\Psi,\varphi) \d t=0, 
\qquad\qquad\qquad\forall\, \varphi\in L^2(0,T;{\cal H}^1).
\end{align}

For any given ${\bf a}\in L^2(0,T;{\bf H}({\rm curl},{\rm div}))$
and $\chi\in L^2(0,T;H^1)$,
we let ${\bf a}_{h,\tau}\in L^2(0,T;{\mathbb N}_h^k)$ and 
$\chi_{h,\tau}\in L^2(0,T;{\mathbb V}_h^{k+1})$ be 
finite element functions such that 
\begin{align*}
&{\bf a}_{h,\tau}\rightarrow {\bf a}
&&\mbox{strongly
in $L^2(0,T;{\bf H}({\rm curl}))$ as $h\rightarrow 0$} ,\\
&\chi_{h,\tau}\rightarrow \chi
&&\mbox{strongly
in $L^2(0,T;H^1)$ as $h\rightarrow 0$} .
\end{align*} 
The equations \eqref{FEM3}-\eqref{FEM2} imply 
\begin{align*}%\label{FFFFF1}
&\int_0^T\bigg[ (\phi_{h,\tau}^+, \chi_{h,\tau}) 
- ({\bf A}_{h,\tau}^+,\nabla \chi_{h,\tau} ) \bigg]\d t= 0 \, , \\[10pt]
&\int_0^T\bigg[(\partial_t{\bf A}_{h,\tau},{\bf a}_{h,\tau}) 
+ (\nabla\phi_{h,\tau}^+ \, ,{\bf a}_{h,\tau})
+ (\nabla\times{\bf A}_{h,\tau}^+  \, , \nabla\times {\bf a}_{h,\tau})\bigg]\d t\nn\\
&\quad +\int_0^T\bigg[ {\rm Re}\bigg(\overline\psi_{h,\tau}^- \bigg(\frac{i}{\kappa} \nabla 
+{\bf A}_{h,\tau}^-\bigg) \psi_{h,\tau}^-  , {\bf a}_{h,\tau}\bigg)\bigg]\d t
=\int_0^T\bigg[ (\nabla\times{\bf H} \, ,{\bf a}_{h,\tau})\bigg]\d t  \, . 
%\label{FFFFF2}
\end{align*}
Let $h=h_m\rightarrow 0$
and $\tau=\tau_m\rightarrow 0$ in 
the last two equations and 
use \eqref{ConvA1} and \eqref{Convpsi1-3}-\eqref{WWWA}.
We obtain 
\begin{align}
&\int_0^T\Big[ (\Phi, \chi) 
-  ({\bf \Lambda} ,\nabla \chi) \Big] \d t = 0 \, , 
\label{PDEPhiF}\\[10pt]
&\int_0^T\bigg[(\partial_t{\bf \Lambda},{\bf a}) 
+ (\nabla\Phi \, , {\bf a})
+ (\nabla\times{\bf \Lambda}  \, , \nabla\times {\bf a})
+ {\rm Re}\bigg(\overline\Psi\bigg(\frac{i}{\kappa} \nabla 
+{\bf \Lambda}\bigg) \Psi  , {\bf a} \bigg)\bigg]\d t\nn\\
&=\int_0^T (\nabla\times{\bf H} \, ,{\bf a} )  \d t \, ,
\label{PDEAF}
\end{align}
which hold for any given 
${\bf a}\in L^2(0,T;{\bf H}({\rm curl},{\rm div}))$
and $\chi\in L^2(0,T;H^1)$.
Since \eqref{PDEPhiF} implies
$\Phi=-\nabla\cdot{\bf\Lambda}$,
\eqref{PDEAF} can be rewritten as
\begin{align}
&\int_0^T\bigg[(\partial_t{\bf \Lambda},{\bf a})  
+ (\nabla\cdot {\bf \Lambda}\, , \nabla\cdot{\bf a})
+ (\nabla\times{\bf \Lambda}  \, , \nabla\times {\bf a})
+ {\rm Re}\bigg(\overline\Psi\bigg(\frac{i}{\kappa} \nabla 
+{\bf \Lambda}\bigg) \Psi  , {\bf a} \bigg) \bigg]\d t\nn\\
& = \int_0^T(\nabla\times{\bf H} \, ,{\bf a} ) \d t  \, ,
\qquad\qquad\qquad\qquad
\qquad\quad\forall\, {\bf a}\in L^2(0,T;{\bf H}({\rm curl},{\rm div})) .
\label{PDEAF2}
\end{align}
From \eqref{PDEPsiFM} and \eqref{PDEAF2}
we see that $(\Psi,{\bf\Lambda})$
is a weak solution of the PDE problem \eqref{PDE1}-\eqref{PDEini} 
with the regularity 
\begin{align*}
&\Psi\in C([0,T];{\mathcal  L}^2)\cap 
L^\infty(0,T;{\mathcal  H}^1) ,
\quad \partial_t\Psi \in L^2(0,T;{\mathcal  L}^2),
\quad |\Psi|\leq 1~~\mbox{a.e.~in~\,}\Omega\times(0,T),\\
& {\bf \Lambda}\in C([0,T];{\bf L}^2)\cap 
L^{\infty}(0,T;{\bf H}({\rm curl},{\rm div})) , 
\quad \partial_t{\bf \Lambda}\in L^2(0,T;{\bf L}^2) . 
\end{align*}
Since the PDE problem \eqref{PDE1}-\eqref{PDE2} 
has a unique weak solution
with the regularity above 
(see appendix), 
it follows that 
$\Psi=\psi$, ${\bf\Lambda}={\bf A}$
and $\Phi=\phi$. 

Overall, 
we have proved that  
any sequence 
$(\psi_{h_m,\tau_m}^+,\phi_{h_m,\tau_m}^+, {\bf A}_{h_m,\tau_m}^+)$
with 
$h_m,\tau_m\rightarrow 0$ contains a  
subsequence which converges to 
the unique solution 
$(\psi,\phi, {\bf A})$ of the PDE problem 
\eqref{PDE1}-\eqref{PDEini} 
in the sense of \eqref{Convpsi1-3}-\eqref{ConvA2-3-4}.
This implies that 
$(\psi_{h,\tau}^+,\phi_{h,\tau}^+, {\bf A}_{h,\tau}^+)$
converges to 
$(\psi,\phi, {\bf A})$
as $h,\tau\rightarrow 0$
in the sense of Theorem \ref{MainTHM}. 

The proof of Theorem \ref{MainTHM}
is complete. \endproof

\section{Numerical example}
\setcounter{equation}{0}

\begin{figure}[h]
\centering
%\subfigure[The domain $\Omega$ (shadow part)] 
{\includegraphics[height=1.5in,width=1.7in]{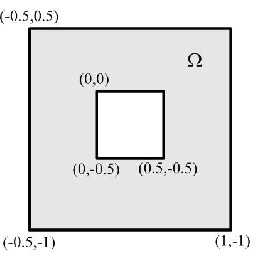}}
\qquad\qquad
%\subfigure[The domain $\Omega_0$ (shadow part)] 
{\includegraphics[height=1.4in,width=2.25in]{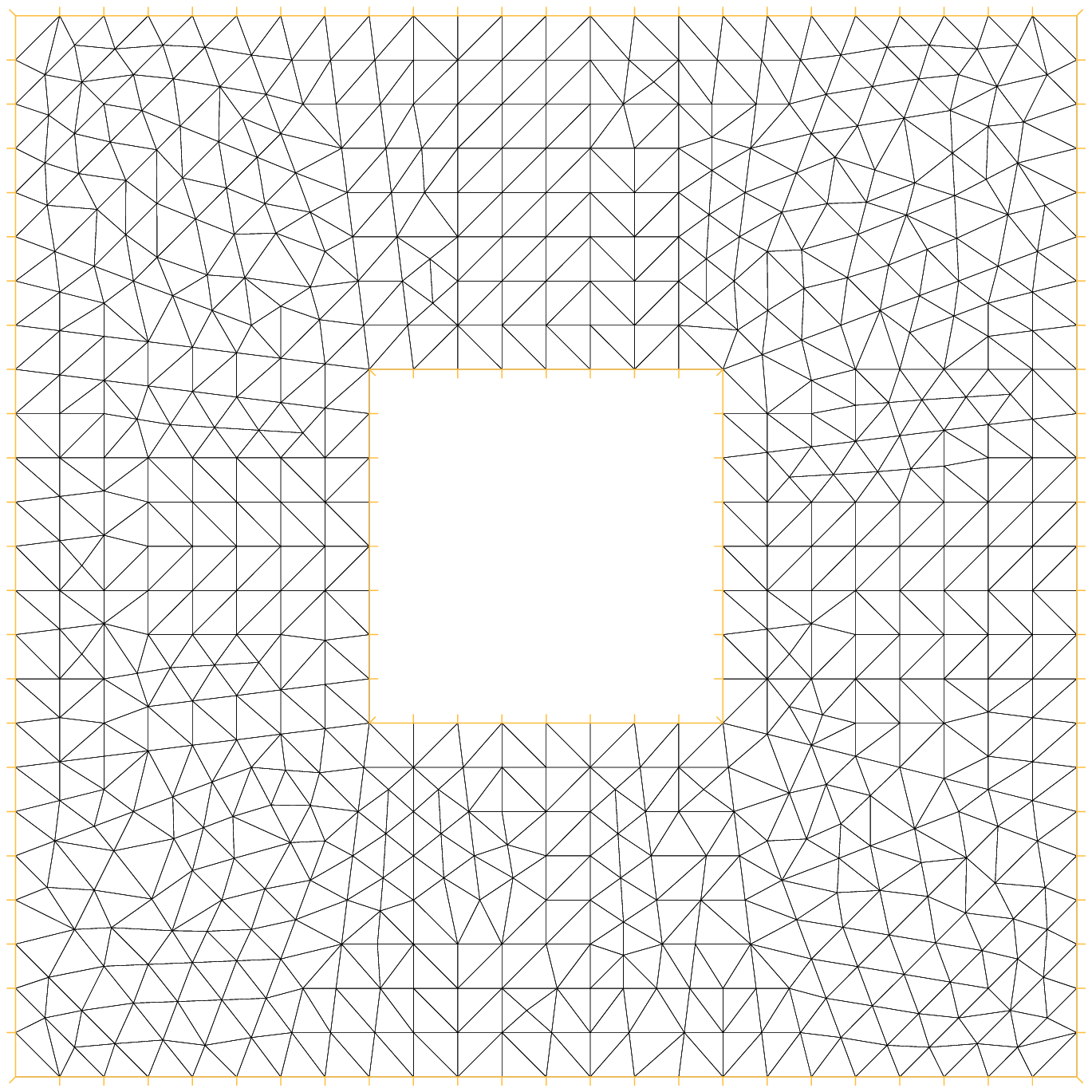}}
\caption{Illustration of the computational domain 
and the triangulation.}
\label{FigD2}
\end{figure}

We consider the equations 
\begin{align}
&\eta\frac{\partial \psi}{\partial t} -i\eta \kappa \psi \nabla\cdot{\bf A}
+ \bigg(\frac{i}{\kappa} \nabla 
+ \mathbf{A}\bigg)^{2} \psi
 + (|\psi|^{2}-1) \psi  = g ,
\label{NTPDE1}\\[5pt]
&\frac{\partial \mathbf{A}}{\partial t} 
-\nabla(\nabla\cdot{\bf A}) 
+ \nabla\times(\nabla\times{\bf A})
+  {\rm Re}\bigg[\overline\psi\bigg(\frac{i}{\kappa} \nabla 
+ \mathbf{A}\bigg) \psi\bigg] 
=  {\bf g}+\nabla\times H ,
\label{NTPDE2}
\end{align}
in a nonsmooth, nonconvex and multi-connected 
two-dimensional domain $\Omega$, as shown in Figure \ref{FigD2}, 
where we use the notations
\begin{align*}
&\nabla\times {\bf A}
=\frac{\partial A_2}{\partial x_1}-\frac{\partial A_1}{\partial x_2},
\qquad 
\nabla\cdot {\bf A}=\frac{\partial A_1}{\partial x_1}
+\frac{\partial A_2}{\partial x_2},\\
&\nabla\times H=\bigg(\frac{\partial H}{\partial x_2},\,
-\frac{\partial H}{\partial x_1}\bigg),\quad 
\nabla\psi=\bigg(\frac{\partial \psi}{\partial x_1},\,
\frac{\partial \psi}{\partial x_2}\bigg).
\end{align*}
The artificial right-hand sides 
$H=\nabla\times{\bf A}\in C([0,T];{\bf H}^2)$, 
$g\in C([0,T];L^2)$ and 
${\bf g}\in C([0,T];{\bf L}^2)$ are chosen 
corresponding to the exact solution 
(written in the polar coordinates)
\begin{align*}
&\psi=t^2\Phi(r)r^{2/3}\cos(2\theta/3),\\
&{\bf A}=\Big (\big( 4t^2 \Phi(r)r^{-1/3}/3
+t^2\Phi'(r)r^{2/3}\big)\cos(\theta/3),~
\big( 4t^2 \Phi(r)r^{-1/3}/3
+t^2\Phi'(r)r^{2/3}\big)\sin(\theta/3)\Big ) ,  
\end{align*}

\begin{table}[h]
\begin{center}
\caption{\small Errors of the Galerkin finite 
element solution 
with $\tau=2h$.} 
\label{Tab1} \medskip
\begin{tabular}{c|c|c|c|cc}
\hline
 $h$ &  $\| \psi_h^N -\psi^N \|_{L^2}$ 
 &  $\| |\psi_h^N| -|\psi^N| \|_{L^2}$ &
$\| {\bf A}_h^N - {\bf A}^N \|_{L^2}$ 
& $\| {\bf B}_h^N - {\bf B}^N \|_{L^2}$        \\
\hline
  1/32
&3.3872E-03  &2.5568E-03  &9.2707E-02 &2.5726E-01  \\
  1/64
&2.9051E-03  &1.7546E-03 &9.1339E-02 &1.7235E-01\\
 1/128
&2.7352E-03  &1.4476E-03 &9.0496E-02 &1.4259E-01 \\
% 1/256
%&2.6504E-03  &1.3098E-03 &9.0010E-02 &1.3423E-01 \\
\hline  
convergence rate &   $O(h^{0.09})$  
& $O(h^{0.29})$ 
& $O(h^{0.01})$  
& $O(h^{0.27})$\\ 
\hline 
\end{tabular}
\end{center}

\begin{center}
\caption{\small Errors of the mixed finite 
element solution with $\tau=2h$.}
\label{Tab2} \medskip
\begin{tabular}{c|c|c|c|cc}
\hline
 $h$ &  $\| \psi_h^N -\psi^N \|_{L^2}$ 
 &  $\| |\psi_h^N| -|\psi^N| \|_{L^2}$  &
$\|  {\bf A}_h^N - {\bf A}^N \|_{L^2}$ 
&$\|{\bf B}_h^N - {\bf B}^N \|_{L^2}$     \\
\hline
  1/32
&5.0142E-03  &2.9762E-03 &4.1846E-03 &  1.7284E-01 \\
  1/64
&1.8455E-03  &1.4828E-03 &2.3881E-03 & 8.7132E-02\\ 
 1/128
&7.5068E-04  &5.6680E-04 &1.4964E-03 &4.3196E-02\\
% 1/256
%& E-04  & E-04 & E-03 & E-02\\
\hline
convergence rate  & $O(h^{1.29})$   
&  $O(h^{1.38})$ 
&  $O(h^{0.67})$  
&  $O(h^{1.01})$  \\
\hline 
\end{tabular}

\end{center}
\end{table}

\noindent where the cut-off function $\Phi(r)$ is given by
$$
\Phi(r)=\left\{
\begin{array}{ll}
0.1 & \mbox{if}~~r<0.1,  \\
\Upsilon(r) &\mbox{if}~~ 0.1\leq r\leq 0.4 ,\\
0 & \mbox{if}~~r>0.4, 
\end{array}\right.
$$
and $\Upsilon(r)$ is the unique $7^{\rm th}$ 
order polynomial satisfying the 
conditions $\Upsilon'(0.1)=\Upsilon''(0.1)
=\Upsilon'''(0.1)=\Upsilon(0.4)=\Upsilon'(0.4)
=\Upsilon''(0.4)=\Upsilon'''(0.4)=0$ 
and $\Upsilon(0.1)=0.1$.
%The exact solution $(\psi,{\bf A})$
%satisfies the boundary and initial conditions
%\eqref{bc}-\eqref{init} with $\psi_0=0$ and ${\bf A}_0=(0,0)$.

%\begin{figure}[h]
%\vspace{0.1in}
%\centering
%\begin{tabular}{ccc}
%\epsfig{file=Th16.eps,height=1.3in,width=1.9in}
%\hspace{-0.2in}
%\epsfig{file=Th32.eps,height=1.3in,width=1.9in}
%\hspace{-0.2in}
%\epsfig{file=Th64.eps,height=1.3in,width=1.9in}
%\end{tabular}
%%\vskip-0.1in
%\caption{\small Quasi-uniform triangulations with $M=16,32,64$.}
%\label{Lshape}
%\end{figure}

We solve \eqref{NTPDE1}-\eqref{NTPDE2}  
by the linear Galerkin FEM 
and our mixed FEM with $r=k=1$, respectively, with 
the same time-stepping scheme
under the same quasi-uniform mesh, 
and present the errors of the numerical solutions
in Table \ref{Tab1}--\ref{Tab2}, 
where $h$ denotes the distance between the mesh nodes on
$\partial\Omega$ and 
the convergence rate of $\psi_h^N$ is calculated 
%by 
%\begin{align*}
%&{\rm convergence~ rate~ of~}\psi_h^N=
%\log(\|\psi_h^N -\psi^N \|_{L^2}/\|\psi_{h/2}^N -\psi^N \|_{L^2})/\log 2 
%\end{align*}
based on the finest mesh size $h$.
We see that the numerical solution of 
the Galerkin FEM does not decrease to zero,
while the mixed finite element solution proposed in
this paper has an explicit convergence rate $O(h^{0.67})$,  
which is consistent with the regularity 
${\bf A}\in L^\infty(0,T;{\bf H}({\rm curl, div}))
\hookrightarrow L^\infty(0,T;{\bf H}^{2/3-\epsilon})$ 
(though we have not proved such explicit 
convergence rate in this paper). \bigskip

\section*{Appendix: Well-posedness of the PDE
problem {\bf(\ref{PDE1})-(\ref{PDEini})}}
\renewcommand{\thelemma}{A.\arabic{lemma}}
\renewcommand{\thetheorem}{A.\arabic{theorem}}
\renewcommand{\theequation}{A.\arabic{equation}} 
\renewcommand{\theremark}{A.\arabic{remark}} 
\setcounter{equation}{0}
\setcounter{theorem}{0} 
\setcounter{remark}{0}

\begin{theorem} 
{\it 
There exists a unique weak solution 
of (\ref{PDE1})-(\ref{PDEini}) with the following regularity: 
\begin{align*} 
&\psi\in C([0,T];{\mathcal  L}^2)\cap 
L^\infty(0,T;{\mathcal  H}^1) ,
\quad \partial_t\psi \in L^2(0,T;{\mathcal  L}^2),
\quad |\psi|\leq 1~~\mbox{a.e.~in~\,}\Omega\times(0,T),\\
& {\bf A}\in C([0,T];{\bf L}^2)\cap 
L^{\infty}(0,T;{\bf H}({\rm curl},{\rm div})) , 
\quad \partial_t{\bf A}\in L^2(0,T;{\bf L}^2) . 
\end{align*} 
} 
\end{theorem} 
\noindent{\it Proof.}$\,\,\,$
From \eqref{PDEPsiFM} and \eqref{PDEAF2}
we see that there exists a weak solution 
$(\Psi,{\bf\Lambda})$ 
of  \eqref{PDE1}-\eqref{PDEini}
with the regularity above. It remains to prove
the uniqueness of the weak solution. 

Suppose that there are two weak solutions 
$(\psi,{\bf A})$ and $(\Psi,{\bf\Lambda})$ for 
the system \eqref{PDE1}-\eqref{PDEini}.
Then we define 
$e=\psi-\Psi$ and ${\bf E}={\bf A}-{\bf\Lambda}$
and consider the difference equations
\begin{align}
&\int_0^T\Big[\big(\eta\partial_t e  ,\varphi\big)
+ \frac{1}{\kappa^2}\big(\nabla  e, \nabla\varphi\big) 
+ \big(|{\bf A}|^2   e,  \varphi\big) \Big]\d t\nn\\
& =\int_0^T\Big[-\frac{i}{\kappa}\big({\bf A}\cdot\nabla e ,\varphi\big)
-\frac{i}{\kappa}\big({\bf E}\cdot\nabla \Psi ,\varphi\big)
+\frac{i}{\kappa}\big(e {\bf A},\nabla\varphi\big)
+\frac{i}{\kappa}\big(\Psi {\bf E},\nabla\varphi\big)  \nn\\
&\quad - \big((|{\bf A}|^2 -|{\bf \Lambda}|^2)  \Psi,  \varphi\big)
 -\big( (|\psi|^{2}-1) \psi-(|\Psi|^{2}-1) \Psi,\varphi\big)\Big]\d t\nn\\
&\quad -\int_0^T\big(i\eta\kappa\psi\nabla\cdot{\bf E}
+i\eta\kappa e\nabla\cdot{\bf \Lambda},\varphi\big)\d t  ,
\label{UErEq1}
\end{align}
and
\begin{align}
&\int_0^T\Big[\big(\partial_t{\bf E} ,{\bf a}\big)
+ \big(\nabla\times{\bf E},\nabla\times{\bf a}\big)
+\big(\nabla\cdot{\bf E},\nabla\cdot{\bf a}\big) \Big]\d t\nn\\
& =-\int_0^T{\rm Re} \bigg( 
\frac{i}{\kappa}( \overline\psi\nabla  \psi- \overline\Psi \nabla  \Psi)
+  {\bf A}(|\psi|^2-|\Psi|^2)+|\Psi|^2 {\bf E}\, ,\, {\bf a}\bigg) \d t ,
\label{UErEq2}
\end{align}
which hold for any $\varphi\in L^2(0,T;{\mathcal H}^1)$ and 
${\bf a}\in L^2(0,T;{\bf H}({\rm curl},{\rm div}))$.
Choosing $\varphi(x,t)=e(x,t)1_{(0,t')}(t)$ in \eqref{UErEq1} 
and considering the real part,
%, by using the inequality
%$$
%{\rm Re}\big(  |\psi|^{2} \psi- |\Psi|^{2}  \Psi,\psi-\Psi\big)\geq 0 ,
%$$
we obtain
\begin{align*}
& \frac{\eta}{2} \|e(\cdot,t') \|_{{\cal L}^2}^2 
+ \int_0^{t'}\Big(\frac{1}{\kappa^2}\|\nabla  e\|_{{\cal L}^2}^2
+  \|{\bf A}    e\|_{{\bf L}^2}^2\Big)\d t \\
&\leq 
\int_0^{t'}\Big(C\|{\bf A}\|_{{\bf L}^{3+\delta}}\|\nabla e\|_{{\cal L}^2}
 \|e\|_{{\cal L}^{6-4\delta/(1+\delta)}}
+C\|{\bf E}\|_{{\bf L}^{3+\delta}}\|\nabla \Psi\|_{{\cal L}^2}
\|e\|_{{\cal L}^{6-4\delta/(1+\delta)}}\\
&\quad 
+C\| e\|_{{\cal L}^{6-4\delta/(1+\delta)}} \|{\bf A}\|_{{\bf L}^{3+\delta}}
\|\nabla e\|_{{\cal L}^2}  +C\| {\bf E}\|_{{\bf L}^2}\|\nabla e\|_{{\cal L}^2}  \\
&\quad
+C(\|{\bf A}\|_{{\bf L}^{3+\delta}}+\|{\bf \Lambda}\|_{{\bf L}^{3+\delta}})
\|{\bf E}\|_{{\bf L}^2} \|e\|_{{\cal L}^{6-4\delta/(1+\delta)}} +C\| e\|_{{\cal L}^2}^2
+C\|\nabla\cdot{\bf E}\|_{L^2}\|e\|_{{\cal L}^2}\Big)\d t \\
&\leq \int_0^{t'}\Big(C\|\nabla e\|_{L^2}
 (C_\epsilon \|e\|_{{\cal L}^2}+\epsilon\|\nabla e\|_{{\cal L}^2})
 +C\|{\bf E}\|_{{\bf H}({\rm curl},{\rm div})}
 (C_\epsilon \|e\|_{{\cal L}^2}+\epsilon\|\nabla e\|_{{\cal L}^2})\\
&\quad
+C\|\nabla e\|_{{\cal L}^2}(C_\epsilon \|e\|_{{\cal L}^2}
+\epsilon\|\nabla e\|_{{\cal L}^2}) 
+C\| {\bf E}\|_{{\bf L}^2}\|\nabla e\|_{{\cal L}^2} \\
&\quad
+C\|{\bf E}\|_{{\bf L}^2}(C_\epsilon \|e\|_{{\cal L}^2}
+\epsilon\|\nabla e\|_{{\cal L}^2}) +C\| e\|_{{\cal L}^2}^2
+C\|\nabla\cdot {\bf E}\|_{L^2}\| e\|_{{\cal L}^2} \Big)\d t\\  
&\leq \int_0^{t'}\Big(\epsilon\|\nabla e\|_{{\cal L}^2}^2+
\epsilon\|\nabla\times{\bf E}\|_{{\bf L}^2}^2
+ \epsilon\|\nabla\cdot{\bf E}\|_{L^2}^2 
+C_\epsilon\|e\|_{{\cal L}^2}^2
+ C_\epsilon\|{\bf E}\|_{{\bf L}^2}^2\Big)\d t ,
\end{align*}
where $\epsilon$ can be arbitrarily small.
By choosing ${\bf a}(x,t)={\bf E}(x,t)1_{(0,t')}(t)$ in 
\eqref{UErEq2}, we get
\begin{align*}
& \frac{1}{2}\|{\bf E}(\cdot,t')\|_{{\bf L}^2}^2 
+\int_0^{t'}\Big(\|\nabla\times{\bf E} \|_{{\bf L}^2}^2
+\|\nabla\cdot{\bf E} \|_{{\bf L}^2}^2 \Big)\d t\\
&\leq \int_0^{t'}\Big(C \| e\|_{{\cal L}^{6-4\delta/(1+\delta)}}
\|\nabla  \psi\|_{L^2}\| {\bf E}\|_{{\bf L}^{3+\delta}}
+C\|\nabla  e \|_{{\cal L}^2} \| {\bf E}\|_{{\bf L}^2}  \\
&\quad
+ (\|e\|_{L^{6-4\delta/(1+\delta)}}\| {\bf A}\|_{{\bf L}^{3+\delta}}
+\|{\bf E}\|_{{\bf L}^2})\|{\bf E}\|_{L^2}\Big)\d t\\
&\leq \int_0^{t'}\Big(C(C_\epsilon \| e\|_{{\cal L}^2} 
+\epsilon\|\nabla  e \|_{{\cal L}^2})
\| {\bf E}\|_{{\bf H}({\rm curl},{\rm div})}
+\|\nabla  e \|_{{\cal L}^2} \| {\bf E}\|_{L^2} \\
&\quad + (\|e\|_{{\cal L}^2} +\|\nabla e\|_{{\cal L}^2}
+\|{\bf E}\|_{{\bf L}^2})\|{\bf E}\|_{{\bf L}^2}\Big)\d t\\
&\leq
\int_0^{t'}\Big(\epsilon\|\nabla e\|_{{\cal L}^2}^2
+\epsilon\|\nabla\times{\bf E}\|_{{\bf L}^2}
+\epsilon\|\nabla\cdot{\bf E}\|_{{\bf L}^2}
+C_\epsilon  \|e\|_{{\cal L}^2}^2
+ C_\epsilon  \|{\bf E}\|_{{\bf L}^2}^2\Big)\d t ,
\end{align*}
where $\epsilon$ can be arbitrarily small.
By choosing $\epsilon<\frac{1}{4}
\min(1, \kappa^{-2} )$ and summing up the  
two inequalities above, we have
\begin{align*}
&  \frac{\eta}{2}\|e(\cdot,t')\|_{L^2}^2
+\frac{1}{2}\|{\bf E}(\cdot,t')\|_{L^2}^2 
\leq 
\int_0^{t'}\Big(C\|e\|_{L^2}^2 
+C\| {\bf E}\|_{L^2}^2\Big)\d t ,
\end{align*}
which implies
\begin{align*}
&\max_{t\in(0,T)}
\bigg(\frac{\eta}{2}\|e\|_{L^2}^2
+\frac{1}{2}\|{\bf E}\|_{L^2}^2\bigg)=0  
\end{align*}
via Gronwall's inequality. 
Uniqueness of the weak solution is proved.
\qed\bigskip

{\bf Acknowledgement.}$\quad$
I would like to express my gratitude 
to Prof. Christian Lubich for the helpful discussions 
on the time discretization,  
and thank Prof. Weiwei Sun for the 
email communications on this topic.  
I also would like to thank Prof. Qiang Du for the communications in 
CSRC, Beijing, on the time-independency 
of the external magnetic field 
and the incompatibility of the initial data 
with the boundary conditions. 
%Finally, I want to thank Dr. Chaoxia Yang 
%for some corrections of the proof and typos. 
%Part of the research was carried out in
%Nanjing University supported by
%the NSFC (grant no. 11301262), and the rest part 
%of the research was carried out in
%Universit\"at T\"ubingen supported by
%the Alexander von Humboldt Foundation. 

 \end{document}